\documentclass[12pt,titlepage]{article}
\usepackage{amsthm,amsfonts,amsmath}
\usepackage{amscd,amssymb,latexsym,epsfig}

\title{The Morse-Witten complex \\
       via dynamical systems
       }

\author{Joa Weber\\
        Universit\"at M\"unchen \\
        Mathematisches Institut\\
        Theresienstr.~39\\
        D-80333 M\"unchen
        }

\date{21 November 2004\\
      Revised: 26 October 2005
}

\newtheorem{theorem}{Theorem}[section]
\newtheorem{corollary}[theorem]{Corollary}

\newtheorem{lemma}[theorem]{Lemma}
\newtheorem{proposition}[theorem]{Proposition}

\newtheorem{definition}[theorem]{Definition}
\newtheorem{remark}[theorem]{Remark}

\newtheorem{example}[theorem]{Example}

\newcommand{\1}{{{\mathchoice {\rm 1\mskip-4mu l} {\rm 1\mskip-4mu l}
{\rm 1\mskip-4.5mu l} {\rm 1\mskip-5mu l}}}}
\renewcommand{\1}{{{\mathchoice {\rm 1\mskip-4mu l} {\rm 1\mskip-4mu l}
{\rm 1\mskip-4.5mu l} {\rm 1\mskip-5mu l}}}}

\newcommand{\N}{{\mathbb{N}}}

\newcommand{\Q}{{\mathbb{Q}}}
\newcommand{\R}{{\mathbb{R}}}

\newcommand{\Z}{{\mathbb{Z}}}
\newcommand{\Bb}{{\mathcal{B}}}

\newcommand{\Ff}{{\mathcal{F}}}

\newcommand{\Ll}{{\mathcal{L}}}   
\newcommand{\Mm}{{\mathcal{M}}}   

\newcommand{\Oo}{{\mathcal{O}}}

\newcommand{\Vv}{{\mathcal{V}}}

\newcommand{\im}{{\rm im\: }}        
\newcommand{\graph}{{\rm graph }}  

\newcommand{\rank}{{\rm rank }}    
\newcommand{\IND}{{\rm ind}}       
 
\newcommand{\Or}{{\rm Or}}            
\newcommand{\Crit}{{\rm Crit}}        
\newcommand{\Hom}{{\rm Hom}}          
\newcommand{\Ho}{{\rm H}}             
\newcommand{\CM}{{\rm CM}}            
\newcommand{\RP}{\R{\rm P}}           
\newcommand{\HM}{{\rm HM}}            
\renewcommand{\d}{{\rm d}}

\newcommand{\norm}{{\rm norm}}

\newcommand{\eps}{{\varepsilon}}

\newcommand{\Cinf}{C^{\infty}}

\def\NABLA#1{{\mathop{\nabla\kern-.5ex\lower1ex\hbox{$#1$}}}}
\def\Nabla#1{\nabla\kern-.5ex{}_{#1}}
\def\Tabla#1{\Tilde\nabla\kern-.5ex{}_{#1}}
\def\abs#1{\mathopen|#1\mathclose|}   
\def\Abs#1{\left|#1\right|}            
\def\norm#1{\mathopen\|#1\mathclose\|}
\def\Norm#1{\left\|#1\right\|}

\renewcommand{\Tilde}{\widetilde}
\renewcommand{\Hat}{\widehat}

\newcommand{\p}{{\partial}}

\begin{document}

\footnotetext[0]{Phone ++49 89 21804534,
                 Fax ++49 89 21804648,
                 Email joa@math.sunysb.edu
}

\maketitle

\vspace*{8cm}
\noindent
{\bf\small Abstract:}
{\small Given a smooth closed
manifold $M$,
the Morse-Witten complex
associated to
a Morse function $f$
and a Riemannian metric $g$ on $M$
consists of
chain groups generated by the
critical points of $f$ and a boundary
operator counting isolated flow lines
of the negative gradient flow.
Its homology reproduces
singular homology of $M$.
The geometric approach
presented here was developed
in~\cite{We93} and
is based on tools from
hyperbolic dynamical systems.
For instance, we apply the
Grobman-Hartman theorem
and the $\lambda$-lemma
(Inclination Lemma) 
to analyze compactness
and define gluing for the
moduli space of flow lines.
}


\tableofcontents

\section{Introduction}
\label{sec:intro}

Throughout let $M$ be a smooth
closed\footnote{Compact 
and without boundary.}
manifold of finite dimension $n$ and
$f$ a smooth function on $M$.
Assume that all critical points
of $f$ are nondegenerate,
denote by $\Crit_k f$
those of Morse
index $k$ and let $c_k$ be
their total number
(for definitions see
Section~\ref{subsec:critical-points}
).
Denoting the \emph{$k^{th}$
Betti number\footnote{Rank of
$\Ho_k(M;\Z)=$
cardinality of a basis
of its free part.}
of $M$} by
$b_k=b_k(M;\Z)$, the
\emph{strong Morse inequalities}
are given by
\begin{equation}\label{eq:M-ineq}
\begin{split}
     c_k-c_{k-1}+\dots\pm c_0
    &\ge b_k-b_{k-1}+\dots\pm b_0,\quad
     k=0,\dots,n-1, \\
     c_n-c_{n-1}+\dots\pm c_0
    &= b_n-b_{n-1}+\dots\pm b_0.
\end{split}
\end{equation}
Consider the
free abelian groups
$\CM_k:=\Z^{\Crit_k f}$,
for $k=0,\dots,n$.
It is well known
that the strong Morse inequalities
are equivalent to the existence
of boundary homomorphisms
$\p_k:\CM_k\to\CM_{k-1}$
whose homology groups are of
rank $b_k$;
see e.g.
\cite{Th49,Sm60,Mi65,Fr79}.

In 1982 Witten~\cite{Wi82}
brought to light a
geometric realization of such a
boundary operator $\p_k$.
Choosing as auxiliary data
a (generic) Riemannian metric $g$ on $M$,
he looked at the negative gradient
flow associated to $(f,g)$.
Given $x\in\Crit_k f$
and $y\in\Crit_{k-1} f$,
there are
only finitely many so
called \emph{isolated}
flow lines running from $x$ to $y$.
Choosing orientations
of all unstable manifolds
one can associate
a characteristic sign $n_u\in\{\pm1\}$
to every isolated flow line $u$.
Witten defined
the boundary operator
$\p_k$ on $x$
by counting all
isolated flow lines 
with signs emanating from $x$.

To simplify matters
one can ignore the
signs by
taking $\Z_2$-coefficients
and counting modulo two.
Here is a first example.

\begin{example}\label{ex:exp1}
\rm
Consider
the manifold shown in
Figure~\ref{fig:fig-exp1}
(for now ignore the gray arrows
indicating orientations).
\begin{figure}[ht]
  \centering
  \epsfig{figure=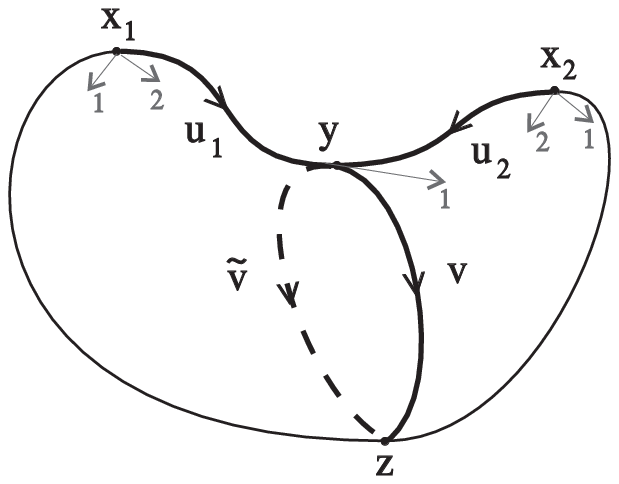}
  \caption{Deformed $S^2$ embedded in $\R^3$}
  \label{fig:fig-exp1}
\end{figure}
The manifold is
supposed to be embedded
in $\R^3$ and the function 
$f$ is given by
measuring height
with respect to the
horizontal coordinate plane.
The metric is induced
by the euclidean metric on
the ambient space $\R^3$.
There are four critical
points $x_1,x_2,y,z$ with
Morse indices $2,2,1,0$,
respectively.
The four isolated flow lines
satisfy
\begin{alignat*}{2}
     \p_2(x_1+x_2)
    &=0\pmod2,
    &
    &
     \\
     \p_1y
    &=0\pmod2,\qquad\quad
    &y
    &=\p_2x_1=\p_2x_2,
    \\
     \p_0z
    &=0,
    &
    &
\end{alignat*}
and the resulting homology
groups
$$
     \HM_2=\langle x_1+x_2 \rangle
     =\Z_2,\qquad
     \HM_1=0,\qquad
     \HM_0=\langle z \rangle
     =\Z_2,
$$
are equal to singular homology
with $\Z_2$-coefficients
of the 2-sphere.
\end{example}

In the early 90's several
approaches towards
rigorously setting up
Witten's complex\footnote{
In~\cite{Wi82}
there is also another
definition via a
deformed deRham complex.
In~1985 Helffer and
Sj\"ostrand~\cite{HS85}
gave a rigorous treatment
using semiclassical
analysis.}
and the resulting
\emph{Morse homology theory}
emerged. The approach by
Floer~\cite{Fl89}
and Salamon~\cite{Sa90}
is via Conley index theory.
The one taken by
Schwarz~\cite{Schw93}
is to consider
the negative gradient
equation in the spirit
of Floer theory as a section
in an appropriate Banach
bundle over the set
of paths in $M$
(see also~\cite{Sa90} for partial
results).
A third approach
from a dynamical systems
point of view, namely
via intersections
of unstable and stable
manifolds, was taken
by the present author.
In~\cite{We93} we
applied the
Grobman-Hartman theorem
and the $\lambda$-lemma
to set up the Morse-Witten
complex. Po\'zniak's
work on the more general
Novikov complex
carries elements of the
second and third approach
and we shall present his
definition of the
continuation maps in
Section~\ref{subsec:continuation}.

Writing the present paper was
motivated by recent
developments.
\emph{Although} unpublished,
the dynamical systems
methods developed
in~\cite{We93}
proved useful in the work
of Ludwig~\cite{Lu03} on
stratified Morse theory.
\emph{Because} 
unpublished,
they were rediscovered
independently by
Jost~\cite{Jo-02}
and -- in the far more general context
of Hilbert manifolds --
by Abbondandolo and
Majer (see~\cite{AM04} and references
therein).

All results
in this paper are part of
mathematical folklore,
unless indicated differently.
In fact,
the proofs in
Sections~\ref{subsec:compactness}--\ref{subsec:gluing}
are due to
the author~\cite{We93}
and the ones in
Section~\ref{subsec:continuation}
to Po\'zniak~\cite{Po91}.

This paper is organized
as follows: Section~\ref{sec:morse-theory}
recalls relevant elements
of Morse theory
and negative gradient flows.
Section~\ref{sec:orbit-spaces}
is at the heart of the matter
where we introduce
and analyze
the moduli spaces
of connecting flow lines.
Morse-Smale transversality
assures that they are manifolds.
Then we show how to
compactify, glue and orient
them consistently
using tools from
dynamical systems
(for which
our main reference will
be the excellent textbook
of Palis and
de~Melo~\cite{PM82}). In
Section~\ref{sec:Morse-homology}
we define
the Morse-Witten complex,
investigate its dependence
on $(f,g)$ and arrive at the
theorem equating its homology
to singular homology.
Finally we provide some
remarks concerning
Morse-inequalities,
relative Morse homology,
Morse cohomology, and
Poincar\'e duality
(for details we
refer the reader
to the original 
source~\cite{Schw93}).
We conclude
by computing Morse
homology and cohomology
of real projective
space $\RP^2$
with coefficients in
$\Z$ and $\Z_2$.

\medskip\noindent
{\small\it Acknowledgements:}
{\footnotesize
At the time of
writing~\cite{We93}
I benefitted from
discussions with
Burkhard Nobbe,
Marcin Po\'zniak,
Dietmar Salamon,
Matthias Schwarz,
S\"onke Seifarth
and Ruedi Seiler.
I am particularly indebted
to Andreas Knauf for
sharing his great expertise
in dynamical systems during
numerous conversations.
Back in the early 90's I recall
warm hospitality
of the Math Departments
at Warwick and Bochum,
in particular of John Jones,
Dietmar Salamon
and Matthias Schwarz.
In present time
I gratefully acknowledge
partial financial
support of
ETH-Projekt Nr.~00321
and of DFG SPP~1154.
I would like to thank
the anonymous referee
for many comments considerably
improving the presentation.
}

\section{Morse theory}
\label{sec:morse-theory}

\subsection{Critical points}
\label{subsec:critical-points}

Given a smooth function
$f:M\to\R$ consider the set
of its critical points
$$
     \Crit f:=
     \{\left. x\in M\right|
     df(x)=0 \}.
$$
Near $x\in \Crit f$ choose
local coordinates
$\varphi=(u^1,\dots,u^n):
U\subset M\to \R^n$
and define a symmetric bilinear
form on $T_xM$, the
\emph{Hessian of $f$ at $x$},
by
$$
     H^f_x(\xi,\eta)
     :=\sum_{i,j=1}^n
     S_{ij}
     \;\xi^i \eta^j,\qquad
     S_{ij}=S_{ij}(f,x;\varphi)
     :=\frac{\p^2 f}{\p u^i \p u^j} (x).
     $$
The symmetric matrix
$S(f,x;\varphi):=
(S_{ij})_{i,j=1}^n$
is the \emph{Hessian matrix}
and the number of its negative eigenvalues
$\IND_f(x)$
is called the \emph{Morse index of $x$}.
If $S(f,x;\varphi)$
is nonsingular, we say that $x$ is
a \emph{nondegenerate} critical point.
It is an exercise to check that
the notions of Hessian,
Morse index and nondegeneracy
do not depend on the choice of
local coordinates
as long as $x\in \Crit f$.
(Hint: Use $df(x)=0$ to show
that for another choice of coordinates
$(\tilde{u}^1,\dots,\tilde{u}^n)$
the matrix $S$ transforms
according to
$\tilde{S}=T^tST$,
where $T$ denotes the
derivative of the coordinate
transition map.
To see that the Morse index
is well defined
apply Sylvester's law; see
e.g.~\cite[Ch.~XV Thm.~4.1]{L95}).

\begin{lemma}\label{le:isolated}
Every nondegenerate critical point
$x$ of $f$ is isolated.
\end{lemma}

\begin{proof}
Choose local coordinates
$(\varphi,U)$
near $x$ as above.
Consider the map
$$
     F
     :=\left(
     \frac{\p f}{\p u^1},\dots,
     \frac{\p f}{\p u^n}\right):
     U\to\R^n
$$
whose zeroes correspond
precisely to the critical
points of $f$ in $U$.
In particular $F(x)=0$.
It remains to show that
there are no other
zeroes nearby.
Since the derivative of
$F$ at $x$
equals $S(f,x;\varphi)$,
it is an isomorphism.
Therefore $F$ is
locally near $x$ a
diffeomorphism
by the inverse function theorem.
\end{proof}

\begin{definition}\rm
A smooth function $f:M\to\R$
is called \emph{Morse}
if all its critical points
are nondegenerate.
\end{definition}

\begin{corollary}\label{co:crit-finite}
If $M$ is closed and
$f$ is a Morse function, then
$\#\Crit f<\infty$.
\end{corollary}

\begin{proof}
Assume not and let
$\{x_k\}_{k\in\N}$
be a sequence of pairwise
distinct critical points of $f$.
By compactness of $M$, there is
a convergent subsequence with limit,
say $x$. By continuity of $df$,
the limit $x$ is again a critical
point. This contradicts
Lemma~\ref{le:isolated}.
\end{proof}

Viewing $df$ as a section of
the cotangent bundle $T^*M$,
the nondegeneracy of a critical
point $x$ is equivalent
to the transversality
of the intersection
of the two closed submanifolds
$M$ and $\graph\: df$ at $x$.
The intersection is compact and,
in the Morse case, also discrete
(complementary dimensions).
This reproves
Corollary~\ref{co:crit-finite}.
Transversality is a generic
property (also open in the case of
closed submanifolds) and so this point
of view is appropriate
to prove the following theorem
(see e.g.~\cite{Hi76}).

\begin{theorem}\label{thm:morse-dense}
If $M$ is closed, then
the set of Morse functions is
open and dense in
$C^\infty(M,\R)$.
\end{theorem}

\subsection{Gradient flows and
(un)stable manifolds}
\label{subsec:grad-flows}

Let $X$ be a smooth vector
field  on $M$. For $q\in M$
consider the initial value problem
for smooth curves
$\gamma:\R\to M$ given by
\begin{equation}\label{eq:grad-flow}
      \dot\gamma(t)
      =X(\gamma(t)), \qquad
      \gamma(0)=q.
\end{equation}
Because $M$ is closed, the solution
$\gamma=\gamma_q$ exists for all
$t\in\R$. It is called the 
\emph{trajectory or flow line
through $q$}.
The \emph{flow generated by $X$}
is the smooth map
$\phi:\R\times M\to M$,
$(t,q)\mapsto\gamma_q(t)$.
For every $t\in\R$, it gives rise to
the diffeomorphism
$\phi_t:M\to M,q\mapsto\phi(t,q)$,
the so called \emph{time-$t$-map}.
The family of time-$t$-maps
satisfies
$\phi_{t+s}=\phi_t \phi_s$
and $\phi_0=id$, i.e. it
is a one-parameter group of
diffeomorphisms of $M$.

The
\emph{orbit $\Oo(q)$ through $q\in M$}
is defined by
$\phi_\R q:=\{\phi_tq\mid t\in\R\}$.
There are three types of orbits,\
namely singular, closed
and regular ones.
A \emph{singular} orbit is one
which consists
of a single point $q$
(which is necessarily a
singularity of $X$).
An orbit is called \emph{closed} 
if there exists
$T\not=0$ such that $\phi_Tq=q$
and $\phi_Tq\not=q$ whenever
$t\in(0,T)$.
In this case $T$ is called the
\emph{period} of the orbit.
Nonsingular and nonclosed orbits are called
\emph{regular}. They are injective immersions
of $\R$ into $M$.
Hence it is natural to ask
if they admit
limit points at their ends.
For $q\in M$, define its $\alpha$-
and \emph{$\omega$-limit} by
\begin{equation*}
\begin{split}
     \alpha(q)
    &:=\{ p\in M\mid\text{
     $\phi_{t_k}q \to p$
     for some sequence
     $t_k\to-\infty$}\}, \\
     \omega(q)
     &:=\{ p\in M\mid\text{
     $\phi_{t_k}q \to p$
     for some sequence
     $t_k\to\infty$}\}.
\end{split}
\end{equation*}
The $\alpha$-limit of $q$
is the $\omega$-limit of $q$
for the vector field $-X$.
Hence the properties of
$\alpha$ translate into
those of $\omega$ and vice versa.
Because $\omega(q)=\omega(\tilde{q})$,
whenever $\tilde{q}$ belongs to the
orbit through $q$, it makes sense to define
$\omega(\Oo(q)):=\omega(q)$.
It is a consequence of
closedness of $M$
that $\alpha(q)$ and $\omega(q)$ are
nonempty, closed, connected
and invariant by the flow
(i.e. a union of orbits);
see e.g.~\cite[Ch.~1 Prop.~1.4]{PM82}.

Let us now restrict to the case of
\emph{gradient flows},
which exhibit a number of
key features.
Let $g$ be a Riemannian
metric and $f$ a smooth function on $M$.
The identity
$g(\nabla f,\cdot)=df(\cdot)$
uniquely determines
the \emph{gradient vector field}
$\nabla f$.
The flow associated to $X=-\nabla f$
is called \emph{negative gradient flow}.
If $\gamma$ is a trajectory of
the negative gradient flow, then
$$
     \frac{d}{dt} f\circ \gamma(t)
     =g(\nabla f(\gamma(t)),\dot\gamma(t))
     =-\Abs{\nabla f(\gamma(t))}^2\le0,\qquad
     \forall t\in\R.
$$
This shows that $f$ is strictly
decreasing along nonsingular orbits.
Therefore closed orbits cannot
exist and any
regular orbit $\Oo(q)$
intersects a level set $f^{-1}(f(q))$
at most once. Moreover,
such an intersection is orthogonal
with respect to $g$.
Using these properties one can show that
$\alpha(q)\cup\omega(q)\subset \Crit f$;
see e.g.~\cite[Ch.~1 Expl.~3]{PM82}.
The idea of proof
is to assume by contradiction
that there exists
$p\in\omega(q)$ with $X(p)\not=0$.
Hence there exists
a sequence $q_k\in\Oo(q)$
converging to  $p$ and 
$f^{-1}(f(p))$ is locally
near $p$ a codimension one submanifold
orthogonal to $X$.
Then, by continuity of the flow,
the orbit through $q$ intersects
the level set in infinitely many points,
which cannot be true.
Example~3 in~\cite[Ch.~1]{PM82}
shows that $\omega(q)$
may indeed contain more than
one critical point.
However, in this case
it must contain infinitely many
by connectedness of $\omega(q)$.
Hence Corollary~\ref{co:crit-finite}
implies Lemma~\ref{le:alpha-lim}
below.

\noindent
The composition of
the linearization of
$-\nabla f$ at a singularity $x$
with the projection onto the second factor
defines the linear operator
$$
     -D\nabla f(x):
     T_xM\stackrel{-d\nabla f(x)}{\longrightarrow}
     T_{\nabla f(x)} TM
     \simeq T_xM\oplus T_xM
     \stackrel{pr_2}{\longrightarrow} T_xM.
$$
With respect to
geodesic normal coordinates
$(\varphi,U)$ near $x$
the operator $D\nabla f(x)$
is represented by the Hessian
matrix
$S(f,x;\varphi)$.
Moreover, these coordinates
are convenient to prove
the identity
$$
     H_x^f(\xi,\eta)
     =g(D\nabla f(x)\:\xi,\eta),\qquad
     \forall \xi,\eta\in T_xM.
$$
Hence $D\nabla f(x)$ is
a symmetric operator
and the number of its negative
eigenvalues coincides with
$k:=\IND_f(x)$.
Let $E^u$ denote the sum
of eigenspaces
corresponding to
negative eigenvalues
and similarly define $E^s$
with respect to positive eigenvalues.
The superscripts abbreviate
\emph{unstable}
and \emph{stable} and this
terminology arises as follows.
The time-$t$-map associated to
the linear vector
field $-D\nabla f(x)$ on
$T_xM$ is given by the
symmetric linear  operator
$A_t:=exp(-tD\nabla f(x))$
on $T_xM$.
Moreover, if $\lambda$ is an eigenvalue
of $D\nabla f(x)$,
then $e^{-t\lambda}$ is eigenvalue of $A_t$
and the eigenspaces are the same.
This shows that $A_t$ leaves the subspaces
$E^u$ and $E^s$ invariant and
acts on them strictly expanding
and contracting, respectively.

\begin{lemma}
For $f\in \Cinf(M,\R)$, let
$\phi_t$ be the time-$t$-map
generated by $X=-\nabla f$.
If $x\in\Crit f$, then
$$
     d\phi_t(x)=exp(-tD\nabla f(x)).
$$
\end{lemma}

\begin{proof}
The map $d\phi_t(x)$
coincides with $A_t$, because it
satisfies the two characterizing
identities for the time-$t$-map
associated to $-D\nabla f(x)$:
pick $\xi\in T_xM$ and
let $c$ be a smooth curve
in $M$ satisfying $c(0)=x$ and
$c'(0)=\xi$. Then use
$\phi_0=id$,
$\p_t\phi_t=-\nabla f(\phi_t)$,
and $\phi_t(x)=x$
to obtain
\begin{equation*}
\begin{split}
     d\phi_0(x)\xi
    &=\left.\frac{\p}{\p\tau}\right|_{\tau=0}
     \phi_0(c(\tau))
     =\left.\frac{\p}{\p\tau}\right|_{\tau=0}
     c(\tau)
     =\xi \\
     \frac{\p}{\p t} d\phi_t(x)\xi
    &=\left.\frac{\p}{\p\tau}\right|_{\tau=0}
     \frac{\p}{\p t}
     \phi_t(c(\tau))
     =-D\nabla f(x) \circ
     d\phi_t(x)\xi.
\end{split}
\end{equation*}
\end{proof}

From now on we shall
assume in addition
that the negative
gradient flow is generated
by a \emph{Morse function}.
Hence, as observed above,
the following lemma
is a consequence of
Corollary~\ref{co:crit-finite}.

\begin{lemma}\label{le:alpha-lim}
Let $M$ be closed and $X=-\nabla f$,
where $f$ is Morse. Then
$\alpha(q)$
and $\omega(q)$
consist each of a single
critical point of $f$,
for every $q\in M$.
\end{lemma}

\begin{figure}[ht]
  \centering
  \epsfig{figure=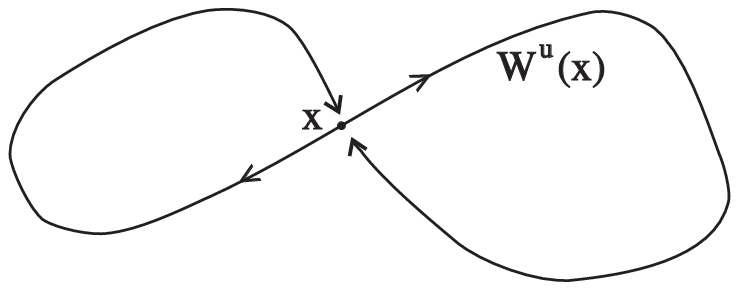}
  \caption{Unstable manifold
           injectively immersed,
           but not embedded}
  \label{fig:fig-fig8}
\end{figure}

The \emph{stable} and
the \emph{unstable manifold
of $x\in\Crit f$} are defined by
$$
     W^s(x):=\{q\in M\mid
     \omega(q)=x\},\qquad
     W^u(x):=\{q\in M\mid
     \alpha(q)=x\}.
$$
Lemma~\ref{le:alpha-lim} shows
$\omega(q)=\lim_{t\to\infty}\phi_tq$.
The map
$
     H:[0,1]\times W^s(x)\to W^s(x)
$,
$
     (\tau,q)\mapsto \phi_{t/(1-t)}q
$,
provides a homotopy between the identity map
on $W^s(x)$
and the constant map $q\mapsto x$.
Hence the (un)stable manifolds
are contractible sets.
Whereas for general
vector fields $X$ with hyperbolic
singularity $x$ these sets are only
injectively immersed (Figure~\ref{fig:fig-fig8}),
they are embedded
in the Morse case.
The reason is, roughly speaking, that
the stable manifold cannot return
to itself, since $f$ is strictly decreasing
along regular orbits.

\begin{theorem}[Stable manifold theorem]
\label{thm:stable-mf}
Let $f$ be a Morse function and
$x\in\Crit f$.
Then $W^s(x)$ is a submanifold
of $M$ without boundary
and its tangent space at $x$ is
given by the stable subspace
$E^s\subset T_xM$.
\end{theorem}

The theorem holds
for $W^u(x)$ with tangent 
space $E^u$
(replace $f$ by $-f$).

\begin{proof}
1) We sketch the proof
that $W^s(x)$ is
locally near $x$
the graph of a smooth map
and $E^s$ is its tangent
space at $x$
(see e.g.~\cite[Ch.~2 Thm.~6.2]{PM82}
for a general hyperbolic flow
and~\cite[Thm.~6.3.1]{Jo-02}
for full details in our case).
Then, using the flow, it follows
that $W^s(x)$ is injectively immersed.
In local coordinates
$(\varphi,U)$
around $x$ the
initial value
problem~(\ref{eq:grad-flow})
is given by
\begin{equation}\label{eq:grad-flow-Rn}
     \dot \gamma
     =A\gamma+h(\gamma),\qquad
     \gamma(0)=\gamma_0.
\end{equation}
Here $\varphi(x)=0$,
the linear map $A$ represents
$-D\nabla f(0)$ in these
local coordinates,
and $h$ satisfies $h(0)=0$
and $dh(0)=0$.
Consider the splitting
$\R^n=E^u\oplus E^s$
induced by $A=(A^u,A^s)$,
fix a metric on $\R^n$
\emph{compatible with
the splitting}\footnote{This means
$\norm{A^s}<1$ and $\norm{A^u}>1$.}
(exists by~\cite[Ch.~2, Prop.~2.10]{PM82}),
and let
$P^u:\R^n\to E^u$
and
$P^s:\R^n\to E^s$
denote the orthogonal projections.
For $\delta,\mu>0$ define
$$
     C_{\delta,\mu}^0
     :=\left\{\gamma\in C^0([0,\infty),\R^n)\left|
     \sup_{t\ge0} e^{\delta t} 
     \Abs{\gamma(t)}\le\mu\right.\right\}.
$$
If $\delta,\mu>0$ are sufficiently small,
then the map
$\Ff:E^s\times 
C_{\delta,\mu}^0\to C_{\delta,\mu}^0$
given by
$$
     \left(\Ff(\gamma_0^s,\gamma)\right)(t)
     :=e^{tA}\gamma_0^s
     +\int_0^t e^{(t-\tau)A}
     P^sh(\gamma(\tau)) \: d\tau
     -\int_t^\infty e^{(t-\tau)A}
     P^uh(\gamma(\tau)) \: d\tau
$$
is a strict contraction
in $\gamma$, whenever $\abs{\gamma_0^s}$
is sufficiently small.
The key fact is that the
unique fixed point
$\hat\gamma$ of
$\Ff(\gamma_0^s,\cdot)$
is precisely the
unique solution
of~(\ref{eq:grad-flow-Rn})
such that
$P^s\hat{\gamma}(0)
=\gamma_0^s$.
Consequently this solution
converges exponentially
fast to zero as $t\to\infty$.
Define the desired graph map
$E^s\to E^u$
locally near zero by
$$
     \gamma_0^s\mapsto
     -\int_0^\infty
     e^{-\tau A}
     P^uh(\hat\gamma(\tau)) \: d\tau.
$$

2) It remains to
prove that $W^s(x)$
is an embedding.
Assume that its
codimension is at least one,
otherwise we are done.
An immersion is locally
an embedding: there exists an open
neighbourhood $W$ of $x$ in $W^s(x)$
which is a submanifold of $M$ of
dimension $\ell:=\dim E^s$.
Let
$
     \mu:=\min_{q\in\p W} f(q)-f(x)
$,
then $\mu>0$ and
$f|_{W^s(x)\setminus W}\ge f(x)+\mu$.
Denote by $B_\eps$ the open
$\eps$-neighbourhood of $x$
with respect to
the Riemannian
distance on $M$.
For $\eps>0$ sufficiently small,
it holds
$f|_{B_\eps}<f(x)+\mu/2$
and therefore
\begin{equation}\label{eq:W1}
     B_\eps\cap (W^s(x)\setminus W)=\emptyset.
\end{equation}
The goal is to construct smooth
submanifold coordinate charts for
every $p\in W^s(x)$.
Assume $p\in W^s(x)\setminus W$,
otherwise we are done.
Define the open neighbourhood
$W_\eps:= W\cap B_\eps$ of $x$ in $W$.
There exists $T>0$ such that
$\phi_Tp\in W_\eps\subset W$.
Choose a submanifold chart $(\varphi,U)$
for $W$ around $\phi_Tp$.
In particular, the set
$U$ is open in $M$, it contains
$\phi_Tp$ and
$\varphi (U\cap W)=0\times V$.
Here $V\subset \R^\ell$ is an open
neighbourhood of $0$.
Shrinking $U$, if necessary, we may assume
without loss of generality that
a) $ U\subset B_\eps$ and
b) $ U\cap W=U\cap W_\eps$.
Condition a)
and~(\ref{eq:W1}) imply
$
     U\cap (W^s(x)\setminus W)=\emptyset
$
and condition b) shows
$$
     U\cap(W\setminus W_\eps)
     =U\cap W \cap (M\setminus W_\eps)
     =U\cap W_\eps \cap (M\setminus W_\eps)
     =\emptyset.
$$
Use these two facts
and represent $W^s(x)$ in the form
$$
     W^s(x)
     =W_\eps
     \cup (W\setminus W_\eps)
     \cup (W^s(x)\setminus W)
$$
to conclude ('no-return')
\begin{equation}\label{eq:Ws}
     U\cap W^s(x)
     =U\cap W_\eps.
\end{equation}
Define the submanifold chart
for $W^s(x)$ at $p$ by
$
     (\psi,U_p)
     :=(\varphi\circ \phi_T,\phi_{-T}U)
$.
The set $U_p$
is indeed an open neighbourhood of $p$
in $M$ and $\psi$ satisfies
\begin{alignat*}{2}
     \psi^{-1}(0\times V)
    &=
     \phi_{-T}\circ\varphi^{-1}
     (0\times V)             
     &&=
     \phi_{-T}(U\cap W)      
    \\
    &=
     \phi_{-T}(U\cap W_\eps) 
    &&=
     \phi_{-T}(U\cap W^s(x)) 
    =U_p\cap W^s(x).
\end{alignat*}
The third equality follows
by condition b)
and equality four by~(\ref{eq:Ws}).
\end{proof}

\section{Spaces of connecting orbits}
\label{sec:orbit-spaces}

Given $x,y\in\Crit f$, define
the \emph{connecting manifold
of $x$ and $y$} by
$$
     \Mm_{xy}
     =\Mm_{xy}(f,g)
     :=W^u(x)\cap W^s(y).
$$
Let $a\in(f(y),f(x))$
be a regular value. The
\emph{space of connecting orbits
from $x$ to $y$} is defined by
\begin{equation}\label{eq:hat-Mm}
     \Hat\Mm_{xy}
     =\Hat\Mm_{xy}(f,g,a)
     :=\Mm_{xy}\cap f^{-1}(a).
\end{equation}
This set represents precisely
the orbits of the negative gradient
flow running from $x$ to  $y$,
because every orbit intersects
the level hypersurface exactly once.
For two different choices of $a$
there is a natural identification
between the corresponding sets
$\Hat\Mm_{xy}(f,g,a)$ which
is provided by the flow.

The structure of this section
is the following.
In Subsection~\ref{subsec:transversality}
we observe that it is possible
to achieve
transversality\footnote{
Two submanifolds $A$ and $B$ of $M$
are said to
\emph{intersect transversally}
if
$$
     T_qA+T_qB=T_qM,\qquad 
     \forall q\in A\cap B.
$$
In this case $A\cap B$
is a submanifold of $M$
whose codimension equals the sum of
the codimensions of $A$ and $B$;
see e.g.~\cite[Ch.1 Thm.~3.3]{Hi76}.}
of all
intersections of
stable and unstable
manifolds simultaneously by
an arbitrarily
small $C^1$-perturbation
of the gradient vector field
within the set of gradient
vector fields.
Then the connecting
manifolds and
the spaces of
connecting orbits are
submanifolds of $M$
without boundary
and their dimensions
are given by
\begin{equation}\label{eq:dim}
     \dim \Mm_{xy}
     =\IND_f(x)-\IND_f(y), \quad
     \dim \Hat\Mm_{xy}
     =\IND_f(x)-\IND_f(y)
     -1.
\end{equation}

In Subsection~\ref{subsec:compactness}
we investigate the structure of the
topological boundary of
the connecting manifolds
and show how this leads to a natural
compactification of
the spaces of connecting orbits.
In case of index
difference $+1$ they
are already compact, hence finite.
The other important
case is index difference $+2$.
Here the connected components of
$\Hat\Mm_{xz}$
are either diffeomorphic
to $S^1$ or to $(0,1)$.
This dichotomy follows from the
fact that these are the only
two types of 1-dimensional
manifolds without boundary.
We shall see that to each end
of an open component
there corresponds
a unique pair of connecting orbits
$(u,v)\in\Hat\Mm_{xy}\times\Hat\Mm_{yz}$,
where $y$ is a critical point of
intermediate index.

The main implication of 
Subsection~\ref{subsec:gluing}
is that every such pair
$(u,v)$ corresponds
to precisely one of the
ends of all open components.
The main ingredient
is the so called gluing map.
More precisely,
we shall define a $C^1$-map
which assigns to $(u,v)$
and to a positive real parameter $\rho$
a unique element of $\Hat\Mm_{xz}$.
Moreover, the limit $\rho\to0$
in the sense of
Section~\ref{subsec:compactness}
corresponds to the original
pair $(u,v)$.

In Subsection~\ref{subsec:orientation}
we prove that a choice of orientations
of all unstable manifolds
induces orientations
of the spaces of connecting orbits
and that they are compatible
with the gluing maps of
Section~\ref{subsec:gluing}.

\subsection{Transversality}
\label{subsec:transversality}

\begin{definition}\label{def:MS}\rm
We say that a gradient vector field
$\nabla^gf$ satisfies
the \emph{Morse-Smale condition}
if $W^u(x)$ and $W^s(y)$ intersect
transversally,
for all $x,y\in\Crit f$.
In this case $(g,f)$ is called a
\emph{Morse-Smale pair}.
\end{definition}

Here is an example which shows
how the Morse-Smale condition
can be achieved by an arbitrarily small
perturbation of the Morse function.

\begin{example}\label{ex:torus-height}
\rm
Consider a 2-torus $T$ embedded
upright in $\R^3$
as indicated in Figure~\ref{fig:fig-tor1}
and let $f:T\to\R$ be given by
measuring height with respect to
the horizontal coordinate plane.
This function admits
four critical points
$M,s_1,s_2$ and $m$ of Morse
indices $2,1,1,0$, respectively.
Let the metric on $T$ be induced from
the ambient euclidean space.
The negative gradient
flow is not Morse-Smale, because 
$W^u(s_1)$ and $W^s(s_2)$
do intersect and therefore
the intersection
cannot be transversal.
However, Morse-Smale
transversality can be achieved
by slightly tilting the torus
as indicated in
Figure~\ref{fig:fig-tor2},
in other words by
perturbing $f$
and thereby destroying 
the annoying flow
lines between $s_1$
and $s_2$.

\begin{figure}[ht]
\begin{minipage}[b]{.46\linewidth}
    \centering
    \epsfig{figure=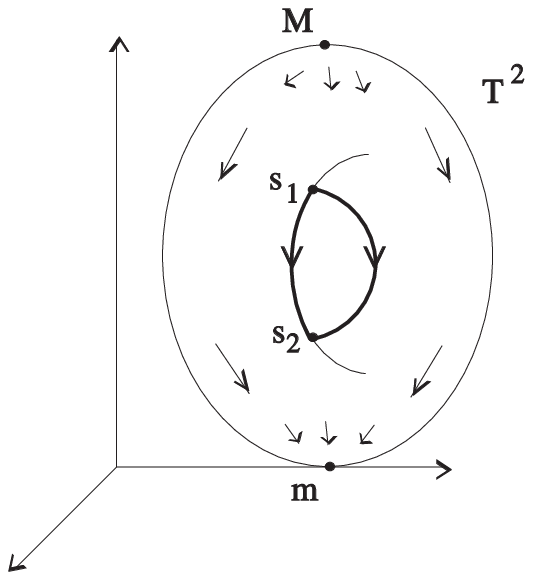}
    \caption{Not Morse-Smale}
    \label{fig:fig-tor1}
\end{minipage}
\hfill
\begin{minipage}[b]{.46\linewidth}
  \centering
  \epsfig{figure=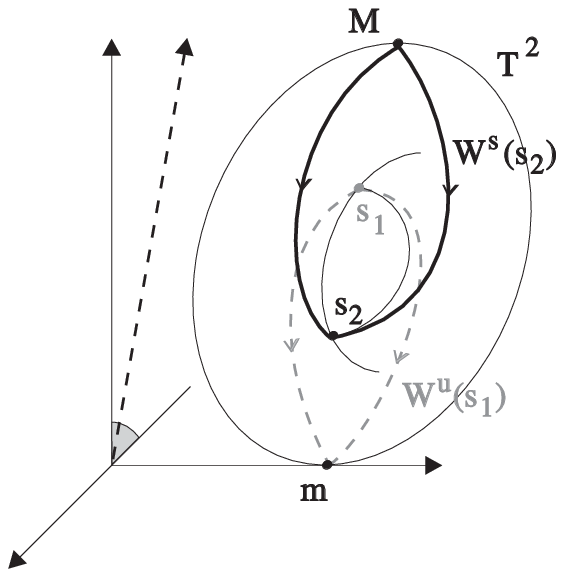}
  \caption{Morse-Smale}
  \label{fig:fig-tor2}
\end{minipage}
\hfill
\end{figure}
\end{example}

\begin{theorem}[Morse-Smale transversality]
\label{thm:MS}
Let  $f$ be
a smooth Morse function and
$g$ a smooth Riemannian
metric on a closed manifold $M$.
Then $\nabla^gf$
can be $C^1$ approximated
by a smooth gradient vector field
$X=\nabla^{\tilde{g}}\tilde{f}$
satisfying the Morse-Smale condition.
\end{theorem}

Theorem~\ref{thm:MS} is due to
Smale~\cite{Sm61}.
Actually $\tilde{f}$
can be chosen such that
its value at any critical
point equals the Morse index.
Note also that the metric
$\tilde{g}$ is generally
not close to $g$ anymore.
It is an exercise to check that
if $(g,f)$ is a Morse-Smale pair,
then $\tilde{f}$ is necessarily Morse.
The \emph{type of a Morse-Smale
vector field $\nabla^gf$},
or equivalently of the
Morse function $f$, is by definition
the number of critical points together
with their Morse indices.
Morse-Smale vector fields
are in particular hyperbolic
vector fields and those have
the property that their type
is locally constant with respect
to the $C^1$-topology.
This can be seen by combining
Proposition~3.1 and the corollary
to Proposition~2.18 in Chapter~2
of \cite{PM82}, which in addition
shows that the critical points
of $f$ and $\tilde{f}$ in
Theorem~\ref{thm:MS} are
$C^0$-close to each other.
In fact, one can even keep the Morse
function $f$ fixed and
achieve Morse-Smale
transversality by perturbing
only the metric
(see~\cite{Schw93}).

\begin{theorem}\label{thm:manifold}
If $-\nabla f$ is Morse-Smale, then
all spaces $\Mm_{xy}$ and
$\Hat\Mm_{xy}$ are submanifolds of $M$
without boundary and their dimensions
are given by~(\ref{eq:dim}).
\end{theorem}

\begin{lemma}\label{le:Mxy}
If $W^u(x)$ and $W^s(y)$
intersect transversally, then
the following are true.
\begin{enumerate}
\item
     If
     $\IND_f(x)<\IND_f(y)$,
     then $\Mm_{xy}=\emptyset$.
\item
     $\Mm_{xx}=\{x\}$.
\item
     If $\IND_f(x)=\IND_f(y)$
     and $x\not=y$, then
     $\Mm_{xy}=\emptyset$.
\item
     If $\Mm_{xy}\not=\emptyset$
     and $x\not=y$, then
     $\IND_f(x)>\IND_f(y)$.
\end{enumerate}
\end{lemma}

\begin{proof}
1. Transversality.
2. Any additional element must be
noncritical and therefore gives
rise to a closed orbit,
which is impossible for gradient flows.
3. Assume the contrary, then
$\Mm_{xy}$ contains
a 1-dimensional
submanifold of the form $\phi_\R q$,
but $\dim\Mm_{xy}=0$.
4. Assume the contrary
and apply statements one and 
three of the lemma
to obtain a contradiction.
\end{proof}

\subsection{Compactness}
\label{subsec:compactness}

Assume throughout this section
that $-\nabla f$ is Morse-Smale
and $x,y\in\Crit f$.
In case that a connecting manifold
$\Mm_{xy}$ is noncompact
we shall investigate the
structure of its topological
boundary as a subset of $M$.
This gives rise to
a canonical compactification
of the associated orbit space
$\Hat\Mm_{xy}$.
In case of index
difference +1, the manifold
$\Hat\Mm_{xy}$ 
itself is already
compact, hence a finite set.
For self-indexing $f$
this is easy to prove.

\begin{proposition}\label{pr:0D-compact}
If $\IND_f(x)-\IND_f(y)=1$,
then
$\#\Hat\Mm_{xy}<\infty$.
\end{proposition}

\begin{proof}
Assume that there is
no critical value between
$f(y)$ and $f(x)$.
If $\Mm_{xy}\not=\emptyset$,
fix $a\in(f(y),f(x))$ and let
$\Hat\Mm_{xy}:=\Mm_{xy}\cap f^{-1}(a)$.
For $\eps>0$ sufficiently small
define two closed sets
$$
     S^u:=f^{-1}(f(x)-\eps)\cap W^u(x),\qquad
     S^s:=f^{-1}(f(y)+\eps)\cap W^s(x).
$$
Let them flow sufficiently long time,
such that $f|_{\phi_TS^u}<a$
and $f|_{\phi_{-T}S^s}>a$
(here we use our assumption).
Being the intersection
of three closed sets, it follows that
the set
$$
     \phi_{[0,T]}S^u \cap
     f^{-1}(a) \cap
     \phi_{[-T,0]}S^s
$$
is closed. On the other hand,
it coincides with the 0-dimensional
submanifold $\Hat\Mm_{xy}$.
A discrete closed
subset of a compact set is finite.
The general case follows
from Theorem~\ref{thm:compactness}
below.
\end{proof}

\begin{definition}\label{def:convergence}\rm
A subset
$K\subset\Hat\Mm_{xy}$ is called
\emph{compact up to broken orbits}, if
\begin{equation}\label{eq:C0-conv}
\begin{split}
    &\text{$\forall$ sequence
     $\{p_k\}_{k\in\N}
     \subset K$}, \\
    &\text{$\exists$ critical points
     $x=x_0,x_1,\dots,x_\ell=y$}, \\
    &\text{$\exists$ connecting orbits
     $u^j\in\Hat\Mm_{x_{j-1}x_j}$,
     $j=1,\dots,\ell$}, \\
    &\text{such that
     $p_k\longrightarrow(u^1,\dots,u^\ell)$ as
     $k\to\infty$}.
\end{split}
\end{equation}
Here convergence means,
by definition,
geometric convergence
with respect to the
Riemannian distance $d$ on $M$
of the orbits through $p_k$
to the union of orbits through
the $u^j$'s.
More precisely,
$$
     \forall \eps>0,
     \exists k_0\in\N,
     \forall k\ge k_0:\:
     \Oo(p_k)\subset
     U_\eps\left(\Oo(u^1)\cup\dots\cup
     \Oo(u^\ell)\right).
$$
Here $U_\eps(A)$ denotes
the open $\eps$-neighbourhood
of a subset $A\subset M$.
We say that
the sequence
\emph{$p_k$ converges
to the broken orbit
$(u^1,\dots,u^\ell)$
of order $\ell$}.
\end{definition}
\begin{figure}[ht]
  \centering
  \epsfig{figure=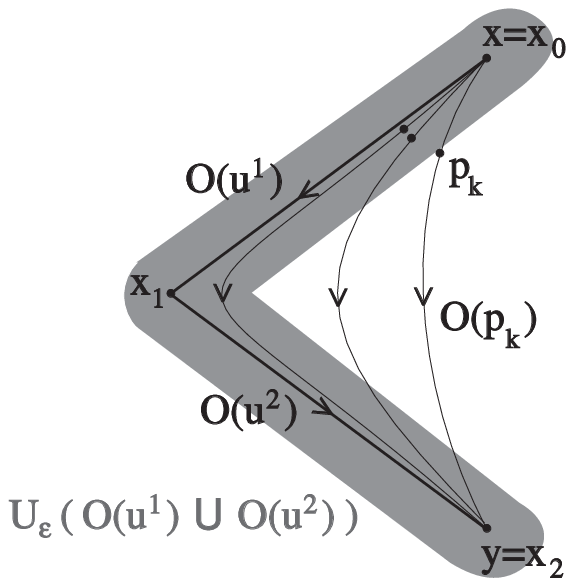}
  \caption{Convergence of
           sequence $p_k$
           to broken orbit
           $(u^1,u^2)$ of
           order 2}
  \label{fig:fig-brok}
\end{figure}
\begin{theorem}[Compactness]
\label{thm:compactness}
If the Morse-Smale condition
is satisfied, then
the spaces of connecting orbits
$\Hat\Mm_{xy}$
are compact up to broken orbits
of order at most
$\IND_f(x)-\IND_f(y)$.
\end{theorem}

\begin{proof}
Fix a regular value $a$ of $f$
and define
$\Hat\Mm_{xy}\subset f^{-1}(a)$
by~(\ref{eq:hat-Mm}).
Assume
$\IND_f(x)>\IND_f(y)$,
otherwise $\Hat\Mm_{xy}=\emptyset$
by Lemma~\ref{le:Mxy} and we are done.
Given a sequence
$\{p_k\}_{k\in\N}\subset
\Hat\Mm_{xy}\subset f^{-1}(a)$,
there exists a subsequence
converging to some element $u$
of the compact set $f^{-1}(a)$.
We use the same notation for the subsequence.
By Lemma~\ref{le:alpha-lim} we have
$u\in\Mm_{z'z}$, for
some $z',z\in\Crit f$.
By continuity of $\phi_t$,
it follows that
$\phi_tu$ lies in the
closure $cl(\Mm_{xy})$ of 
$\Mm_{xy}$
for every $t\in\R$, and therefore
$z\in cl(\Mm_{xy})$.
The proof proceeds in
two steps.

\vspace{.15cm}
\noindent
{\sc Step 1.}
{\it If $z\not=y$, then there exists
$v\in W^u(z)\cap cl(\Mm_{xy})$
with $v\not= z$.
}

\vspace{.1cm}
\noindent
The key tool is the
Grobman-Hartman theorem for flows
which states that the flows associated
to $-\nabla f$ and
$-D\nabla f(z)$, respectively,
are locally conjugate.
(See e.g.~\cite[Ch.~2 Thm.~4.10]{PM82}
where only locally equivalent is
stated, but in fact locally conjugate
is proved; see also~\cite[Thm.5.3]{Rb95}).
This means that there exist
neighborhoods $U_z$ of $z$ in $M$
and $V_0$ of $0$ in $T_zM$,
as well as a homeomorphism
$h:U_z\to V_0$, such that
\begin{equation}\label{eq:loc-conj}
     h(\phi_tq)
     =(D\phi_t(z)\circ h)(q),
\end{equation}
for all $(q,t)$
such that $\phi_tq\in U_z$
and $D\phi_t(z)\circ h(q)\in V_0$.
Observe that $h$
identifies a neighbourhood of
$z$ in $W^s(z)$ with one
of zero in $E^s$;
similarly for the unstable spaces.
(If the eigenvalues of $-D\nabla f(z)$
satisfy certain nonresonance conditions,
then $h$ can be chosen to
be a diffeomorphism; see~\cite{St58}).
We may assume without loss
of generality that $u$ and the $p_k$
are elements of $U_z$, otherwise
apply $\phi_T$ with $T>0$
sufficiently large and choose
a subsequence.
Now apply the Grobman-Hartman
homeomorphism $h$
and consider the image of $u$ and
of the $p_k$ in $V_0\subset T_zM$.
We continue using the same notation
(see Figure~\ref{fig:fig-grha}).
\begin{figure}[ht]
  \centering
  \epsfig{figure=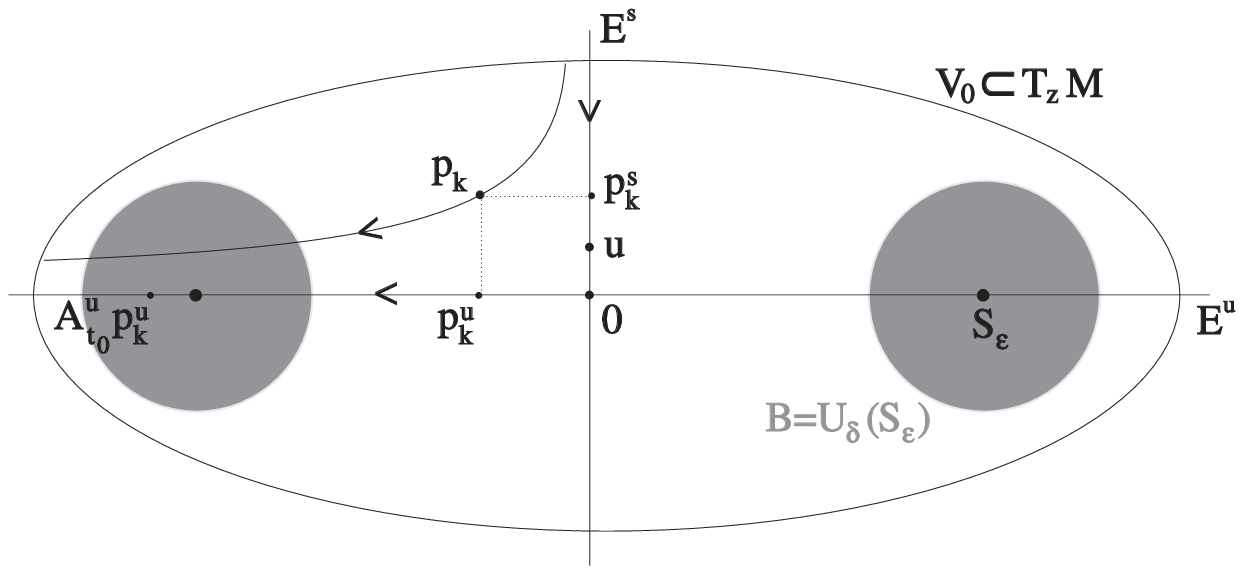}
  \caption{The flow in a Grobman-Hartman
           ($C^0$-) chart}
  \label{fig:fig-grha}
\end{figure}
To prove Step~1 assume the contrary.
Since $h$ conjugates $\phi_t$
and the linearized flow,
the contrary
means that every sphere
$S_\eps$ of radius $\eps$ in $E^u$
admits a $\delta$-neighbourhood
$B$ in $T_zM$ which is
disjoint from $\Mm_{xy}$.
Fix $\eps>0$ and $0<\delta<\eps$
sufficiently small, such that
$S_\eps$ and $B$ are contained in $V_0$.
We may also assume
$\Abs{u}<\delta/2$, otherwise
apply $\phi_T$ again.
By linearity of the flow on
$T_zM=E^s\oplus E^u$
we can write $D\phi_t(z)p_k$
in the form
$(A_t^sp_k^s,A^u_tp_k^u)$,
where $p_k=(p_k^s,p_k^u)$.
The linear operators
$A^s_t\in\Ll(E^s)$
and $A^u_t\in\Ll(E^u)$
introduced in
Section~\ref{subsec:grad-flows}
are, for $t>0$, a strict contraction
and a strict dilatation,
respectively.
For every sufficiently large $k$ we have
$\Abs{p_k^s}<\delta$, hence
$\Abs{A^s_tp_k^s}<\delta$
for positive $t$.
Because $A^u_t$ is expanding
and $0$ is the only fixed point of $A_t$,
it follows that
$\Abs{A^u_{t_0}p_k^u}>\eps$
for some $t_0>0$.
Hence the orbit
$\Oo(p_k)$ runs
through $B$ and this
contradicts our assumption
and therefore proves Step~1.
Furthermore, the
argument shows that
the orbit through $p_k$
converges locally near $z$
to $(u,v)$ in the sense
of~(\ref{eq:C0-conv}).

\vspace{.15cm}
\noindent
{\sc Step 2.}
{\it We prove the theorem.
}

\vspace{.1cm}
\noindent
Assume $z\not=y$, then
by Lemma~\ref{le:alpha-lim}
we conclude that
$v\in\Mm_{z\tilde z}$, for
some $\tilde z\in\Crit f$
with $\IND_f(\tilde z)<\IND_f(z)$.
Repeating the arguments
in the proof of Step~1
leads to an iteration
which can only terminate at $y$.
It must terminate, because
$\Crit f$ is a finite set
by Corollary~\ref{co:crit-finite}
and the index in each
step of the iteration
strictly decreases
by Lemma~\ref{le:Mxy} (d).
Start again with the sequence
$\{p_k\}_{k\in\N}$
and repeat the same arguments
for the flow in backward time.
This proves existence of critical points
and connecting orbits as
in Definition~\ref{def:convergence}.

It remains to prove
uniform convergence.
Near critical points the argument
was given in the proof of Step~1.
Outside fixed neighbourhoods
of the critical points
this is a consequence of the estimate
\[
     d\left(\phi_t q,\phi_t \tilde{q}\right)
     \le e^{\kappa\Abs{t}}
     d\left( q,\tilde{q}\right),\quad
     \forall q,\tilde{q}\in M,
     \forall t\in\R,
\]
where $\kappa=\kappa(M,-\nabla f)>0$ is a constant;
see e.g.~\cite[Ch.~2 Lemma~4.8]{PM82}.
The estimate shows that
on compact time intervals
the orbits through $p_k$ converge
uniformly to the orbit through $u$.
Now set $b:=f^{-1}(v)$
and view $\Hat\Mm_{xy}$
as a subset of $f^{-1}(b)$.
Every point $p_k$ determines
a unique point $\tilde{p}_k$
by intersecting
the orbit through $p_k$
with $f^{-1}(b)$.
Arguing as above, including
choosing further subsequences,
shows that the orbits through the points
$\tilde{p}_k$ converge uniformly
on compact time intervals
to the orbit through $v$.
Repeating this argument
a finite number of times
concludes the proof
of Theorem~\ref{thm:compactness}.
\end{proof}

\subsection{Gluing}
\label{subsec:gluing}

\begin{theorem}[Gluing]\label{thm:gluing}
Assume the Morse-Smale condition is satisfied
and choose $x,y,z\in\Crit f$ of Morse
indices $k+1,k,k-1$, respectively.
Then there exists a positive real number $\rho_0$
and an embedding
$$
     \#:\Hat\Mm_{xy}\times[\rho_0,\infty)\times
     \Hat\Mm_{yz}\to\Hat\Mm_{xz},\qquad
     (u,\rho,v)\mapsto u\#_\rho v,
$$
such that
$$
     u\#_\rho v 
     \longrightarrow (u,v)\quad
     \text{as $\rho\to\infty$}.
$$
Moreover, no sequence in
$\Hat\Mm_{xz}\setminus
(u\#_{[\rho_0,\infty)} v)$
converges to $(u,v)$.
\end{theorem}

\begin{proof}
The proof has three steps.
It consists of local
constructions near $y$.
Therefore we restrict
to the case where
$\phi_t$ is defined near
$y=0\in\R^n$.

\vspace{.15cm}
\noindent
{\sc Step 1 (Local Model).}
{\it We may assume without loss
of generality that a sufficiently
small neighbourhood of $y$
in the stable manifold
is a neighbourhood of $0$ in $E^s$
and similarly for the unstable manifold.
}

\vspace{.15cm}
\noindent
The stable subspace $E^s$
associated to 
$d\phi_t(y)\in\Ll(\R^n)$
is independent of the choice of $t>0$
and similarly for the unstable
subspace $E^u$.
By Theorem~\ref{thm:stable-mf}
they are the tangent spaces at $y$
to the stable and unstable manifold
$W^s$ and $W^u$ of $y$, respectively.
The proof of the theorem
shows that locally
near $y$ the stable and unstable
manifolds are graphs.
More precisely, there exist
small neighbourhoods
$U^s\subset E^s$
and $U^u\subset E^u$ of $y$
and smooth maps
$\eta_s:U^s\to E^u$ and
$\eta_u:U^u\to E^s$ such that
$\eta_s(0)=0$, $d\eta_s(0)=0$
and similarly for $\eta_u$
(see~Figure~\ref{fig:fig-loc1}).
\begin{figure}[ht]
  \centering
  \epsfig{figure=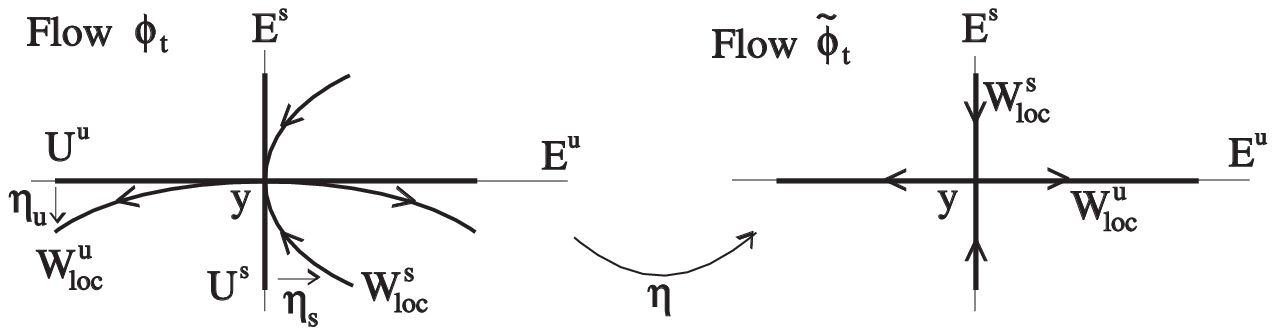,width=\linewidth}
  \caption{Local stable and unstable
           manifolds}
  \label{fig:fig-loc1}
\end{figure}
The graphs of $\eta_u$ and $\eta_s$,
denoted by $W^u_{loc}$
and $W^s_{loc}$,
are called \emph{local} unstable 
and stable manifold, respectively.
The smooth map
$$
     \eta:U^u\times U^s
     \to E^u\oplus E^s,\qquad
     \left(x_u,x_s)\mapsto
     (x_u-\eta_s(x_s),
     x_s-\eta_u(x_u)\right)
$$
satisfies $\eta(0)=0$
and $d\eta(0)=\1$
(see~\cite[Section~7]{PM82}).
Hence it is
a diffeomorphism when restricted to
some neighbourhood of zero.
The family of
local diffeomorphisms
defined by $\tilde{\phi}_t
:=\eta\circ\phi_t\circ\eta^{-1}$, $t>0$,
satisfies $\tilde{\phi}_t(0)=0$
and $d\tilde{\phi}_t(0)=d\phi_t(0)$.
Moreover, a small neighborhood
of zero in the stable manifold of
$\tilde{\phi}_t$
is a small neighbourhood
of zero in $E^s$ and a
similar statement holds
for the unstable manifold.

For later reference
we shall fix a metric $\abs{\cdot}$
on $\R^n$ compatible
with the splitting
$\R^n=E^u\oplus E^s$,
as in the proof of
Theorem~\ref{thm:stable-mf}.

\vspace{.15cm}
\noindent
{\sc Step 2 (Unique intersection point).}
{\it Fix closed balls
$B^u\subset W^u_{loc}\subset E^u$
and $B^s\subset W^s_{loc}\subset E^s$
around $y$ and let
$V:=B^u\times B^s$.
Choose $u\in\Hat\Mm_{xy}$
and assume without loss of generality
that $u\in B^u$ (otherwise replace
$u$ by $\phi_Tu$ for some
$T>0$ sufficiently large).
Choose a $k$-dimensional
disc $D^k\subset W^u(x)$
which transversally
intersects the orbit
through $u$ precisely at $u$.
For $t\ge0$ let $D^k_t$
denote the connected component of
$\phi_t (D^k)\cap V$
containing $\phi_t(u)$.
Choose $v\in\Hat\Mm_{yz}$
and define the $(n-k)$-dimensional
disc $D^{n-k}_{-t}\subset W^s(z)$ similarly,
but with respect to the backward
flow $\phi_{-t}$
(see Figure~\ref{fig:fig-disc}).
Then there exists $t_0\ge0$
such that for every $t\ge t_0$
there is a unique
point $p_t$ of intersection of
$D^k_t$ and $D^{n-k}_{-t}$.
}
\begin{figure}[ht]
  \centering
  \epsfig{figure=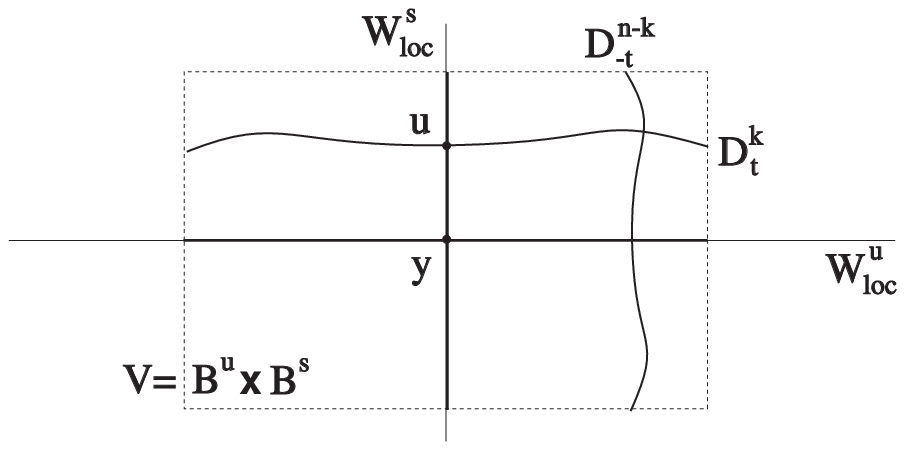}
  \caption{Unique intersection point
           $p_t\in D^k_t\cap D^{n-k}_{-t}
           \subset W^u(x)\cap W^s(z)$}
  \label{fig:fig-disc}
\end{figure}

\vspace{.15cm}
\noindent
The idea of proof
is to represent $D^k_t$
and $D^{n-k}_{-t}$,
for $t>0$ sufficiently large,
as graphs of smooth maps
$F_t:B^u\to B^s$
and $G_t:B^s\to B^u$, respectively.
Since
$D^k_t\cap D^{n-k}_{-t}$
corresponds to the fixed point
set of $G_t\circ F_t:B^u\to B^u$,
it remains to prove that this map is a
strict contraction.

Because $D^k$ intersects $B^s$
transversally, we are in position
to apply our
key tool, namely the
$\lambda$-Lemma~\cite[Ch.~2 Lemma~7.2]{PM82}.
Given $\eps>0$, it asserts
existence of $t_0>0$ such that
$D^k_t$ is $\eps$
$C^1$-close to $B^u$,
for every $t\ge t_0$
(and similarly for $D^{n-k}_{-t}$
and $B^s$).
This means that there exist
diffeomorphisms
\begin{alignat*}{2}
     \varphi_t:B^u
    &\to D^k_t
     ,\qquad
    &q
    &\mapsto(\varphi_t^u(q),\varphi_t^s(q)),\\
     \gamma_t:B^s
    &\to D^{n-k}_{-t}
     ,\qquad
    &p
    &\mapsto(\gamma_t^u(p),\gamma_t^s(p)),
\end{alignat*}
such that
$$
     \Abs{
     \begin{pmatrix}q\\0\end{pmatrix}
     -\begin{pmatrix}\varphi_t^u(q)
     \\\varphi_t^s(q)\end{pmatrix}}
     <\eps,\qquad
     \Norm{
     \begin{pmatrix}\1\\0\end{pmatrix}
     -\begin{pmatrix}d\varphi_t^u(q)
     \\d\varphi_t^s(q)\end{pmatrix}}
     <\eps,\qquad
     \forall q\in B^u,
$$
$$
     \Abs{
     \begin{pmatrix}0\\p\end{pmatrix}
     -\begin{pmatrix}\gamma_t^u(p)
     \\\gamma_t^s(p)\end{pmatrix}}
     <\eps,\qquad
     \Norm{
     \begin{pmatrix}0\\\1\end{pmatrix}
     -\begin{pmatrix}d\gamma_t^u(p)
     \\d\gamma_t^s(p)\end{pmatrix}}
     <\eps,\qquad
     \forall p\in B^s.
$$
Hence, for $\eps>0$ sufficiently
small, the maps $\varphi_t^u$
and $\gamma_t^s$ are invertible
and the required graph maps are given by
(see Figure~\ref{fig:fig-grap})
$$
     F_t(q):=
     \varphi_t^s\circ(\varphi_t^u)^{-1}(q),\qquad
     G_t(p):=
     \gamma_t^u\circ(\gamma_t^s)^{-1}(p).
$$

\begin{figure}[ht]
  \centering
  \epsfig{figure=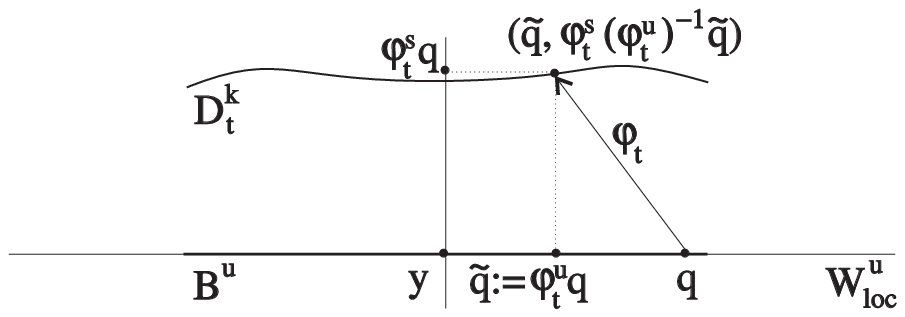}
  \caption{The disc $D^k_t$
           as graph of the map
           $F_t:=
           \varphi_t^s\circ(\varphi_t^u)^{-1}$}
  \label{fig:fig-grap}
\end{figure}

\noindent
The existence of a fixed point
of the smooth map
$G_t\circ F_t:B^u\to B^u$
follows, for instance,
from the Brouwer fixed point
theorem or the fact that
$D_t^k$ and $D_{-t}^{n-k}$
are homotopic to $B^u$ and
$B^s$, respectively, and the
latter have intersection number
one (see e.g.~\cite[Ch.~5 Sec.~2]{Hi76}
for a definition of the intersection
number in case of manifolds
with boundary). Using again
$\eps$ $C^1$-closedness
we obtain
\begin{equation*}
\begin{split}
     \Norm{d\left(G_t\circ 
     F_t\right)|_q}
    &=
     \Norm{
     d\gamma_t^u |_{(\gamma_t^s)^{-1}F_tq}
     \circ d(\gamma_t^s)^{-1}|_{F_tq}
     \circ d\varphi_t^s|_{(\varphi_t^u)^{-1}q}
     \circ d(\varphi_t^u)^{-1}|_q
     } \\
    &\le
     \eps^2 \Norm{
     d(\gamma_t^s)^{-1}|_{F_tq}}
     \cdot \Norm{
     d(\varphi_t^u)^{-1}|_q} \\
    &\le 
     \eps^2/(1-\eps)^2.
\end{split}
\end{equation*}
The last expression is
strictly less than one,
whenever $0<\eps<1/2$.
To obtain the final step
let $S:=d(\varphi_t^u)^{-1}|_q$
and apply the triangle inequality
to obtain
$$
     1
     =\norm{S^{-1}-S^{-1}(\1-S)}
     \ge \norm{S^{-1}}-
     \norm{S^{-1}}\cdot
     \norm{\1-S},
$$
hence an
estimate for $\Norm{S^{-1}}$.
Then the contracting mapping
principle
(see e.g.~\cite[Thm.~V.18]{RS80})
guarantees a \emph{unique}
fixed point and
$\abs{p_t}<\sqrt{2}\eps$.

\vspace{.15cm}
\noindent
{\sc Step 3 (Gluing map).}
{\it Using the notation of
Step~2 we define $\rho_0:=t_0$
and $u\#_\rho v :=p_\rho$,
$\forall \rho\in[\rho_0,\infty)$.
This map satisfies
the assertions of the 
theorem.
}

\vspace{.15cm}
\noindent
The negative gradient
vector field is transverse
to the discs
$D^k_t$ and $D^{n-k}_{-t}$
(otherwise choose $D^k$ and
$D^{n-k}$ in Step~2
smaller). This implies that they
are displaced from themselves
by the flow, so their
intersection $p_t$
cannot remain constant:
$\frac{d}{dt}p_t\not=0$.
This shows that
$u\#_\cdot v$ is an immersion
into $\Mm_{xz}$.
In order to show that it
is an immersion into
$\Hat\Mm_{xz}$ we need to make
sure that $p_t$ does \emph{not}
vary along flow lines, in other words
$\frac{d}{dt}p_t$ and
$-\nabla f(p_t)$ need to be linearly
independent. This is true, since otherwise
the discs must be either
\emph{both} moved in direction
$-\nabla f$ or \emph{both} opposite to it.
However, in our case $D^k_t$ moves in
direction $-\nabla f$
and $D^{n-k}_{-t}$ opposite to it
(see Figure~\ref{fig:fig-embd}).
Hence  $\dim\Hat\Mm_{xz}=1$
implies that the map
$u\#_\cdot v:[t_0,\infty)\to\Hat\Mm_{xz}$
is also a homeomorphism
onto its image (no self-intersections
or returns) and therefore
an embedding.
\begin{figure}[ht]
  \centering
  \epsfig{figure=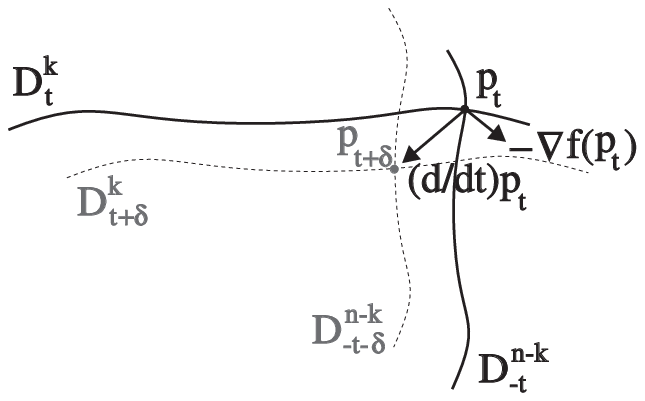}
  \caption{Variation of the point
           of intersection $p_t$}
  \label{fig:fig-embd}
\end{figure}

Choose a sequence of positive
reals $\eps_\ell\to0$.
The $\lambda$-Lemma
yields a sequence
$\{t_{0,\ell}\}_{\ell\in\N}$
such that $D^k_t$ is
$\eps_\ell$ $C^1$-close
to $B^u$, whenever
$t\ge t_{0,\ell}$, and
similarly for $D^{n-k}_{-t}$.
Given any sequence
$t_i\to\infty$ of sufficiently
large reals,
we can choose a subsequence
$\{t_{0,\ell_i}\}_{i\in\N}$
such that
$t_i\ge t_{0,\ell_i}$.
It follows that
$D^k_{t_i}$ and $B^u$
are $\eps_{\ell_i}$
$C^1$-close and similarly
for $D^{n-k}_{-t_i}$ and $B^s$.
Hence
$\abs{p_{t_i}} <\sqrt{2}
\eps_{\ell_i}\to 0$,
as $i\to \infty$.
This proves
$$
     \Abs{p_t}\to 0,\quad
     \text{as $t\to 0$}.
$$
Convergence of $u\#_\cdot v$
to the broken orbit
$(u,v)$ then follows
by the same arguments
as in the proof of
Theorem~\ref{thm:compactness}.
Uniqueness of the intersection
point $p_t$ proves the
final claim of
Theorem~\ref{thm:gluing}.
\end{proof}

\subsection{Orientation}
\label{subsec:orientation}

Assume the Morse-Smale
condition is satisfied.
Recall from
Section~\ref{subsec:grad-flows}
that the stable
and unstable manifolds
are contractible and
therefore orientable.

\begin{proposition}[Induced orientation]
\label{prop:orientation}
Fix an orientation
of $W^u(x)$
for every $x\in\Crit f$
of Morse index larger than zero.
Then, for all
$x,y\in \Crit f$,
the connecting
manifolds $\Mm_{xy}$
and orbits spaces
$\Hat\Mm_{xy}$ \emph{inherit}
induced orientations
$[\Mm_{xy}]_{ind}$ and
$[\Hat\Mm_{xy}]_{ind}$.
\end{proposition}

\begin{proof}
The main idea is that
the transversal intersection
of an oriented and
a cooriented submanifold
is orientable\footnote{
In contrast to the case
of two orientable
submanifolds: consider
$\RP^4=\{[x_0,\dots,x_4]\}$
with submanifolds
$\{[0,x_1,x_2,x_3,x_4]\}\simeq\RP^3$ and
$\{[x_0,x_1,x_2,x_3,0]\}\simeq\RP^3$
intersecting in 
$\{[0,x_1,x_2,x_3,0]\}\simeq\RP^2$.
Here, by definition, an equivalence
class consists of all
vectors which become equal
after scalar multiplication with
some nonzero real number.
}.
An orientation of $W^u(x)$
is by definition an
orientation of its tangent 
bundle\footnote{
A good reference concerning
orientations of vector bundles
is~\cite{Hi76}.}.
Transversality of the intersection
implies that this tangent
bundle splits
along $\Mm_{xy}$. This means
\begin{equation}\label{eq:splitting}
     T_{\Mm_{xy}}W^u(x)
     \simeq
     T\Mm_{xy}\oplus
     \Vv_{\Mm_{xy}}W^s(y).
\end{equation}
Here the last term
denotes the normal
bundle of $W^s(y)$
restricted to $\Mm_{xy}$.
It remains to show
that two of the vector
bundles
are oriented and therefore
determine an orientation
of $T\Mm_{xy}$
denoted by 
$[\Mm_{xy}]_{ind}$.
Firstly, the restriction
of the oriented vector bundle
$TW^u(x)$ to any submanifold,
for instance to
$\Mm_{xy}$,
is an oriented vector bundle.
Secondly, contractibility of the
base manifold $W^s(y)$
implies that $\Vv W^s(y)$
is orientable. Hence an orientation
is determined by an orientation
of a single fiber. The natural choice
is the one over $y$, because
it is isomorphic
to the oriented vector space
$T_y W^u(y)$ via
$$
     T_y W^u(y)\oplus T_y W^s(y)
     \simeq T_yM
     \simeq \Vv_y W^s(y)\oplus T_y W^s(y).
$$
Now restrict the oriented
vector bundle $\Vv W^s(y)\to W^s(y)$
to the submanifold
$\Mm_{xy}\hookrightarrow W^s(y)$
to obtain the required oriented
vector bundle $\Vv_{\Mm_{xy}} W^s(y)$.
The orientation 
$[\Hat\Mm_{xy}]_{ind}$
is determined by the splitting
$$
     T_{\Hat\Mm_{xy}} \Mm_{xy}
     \simeq 
     \R\oplus T\Hat\Mm_{xy},
$$
as the orientation
$[T_{\Hat\Mm_{xy}} \Mm_{xy}]_{ind}$
is given by restriction
of the oriented bundle $T\Mm_{xy}$,
and the orientation of the line bundle
is provided by $-\nabla f$.
\end{proof}

Given $x,y,z\in\Crit f$
of Morse indices $k+1,k,k-1$,
respectively, as well as
$u\in\Hat\Mm_{xy}$ and
$v\in\Hat\Mm_{yz}$,
the gluing map of orbits
$u\#_\rho v=p_\rho$
induces a gluing map
of orientations
$$
     \sigma^\#:
     \Or(\Mm_{xy}^u)\times
     \Or(\Mm_{yz}^v)\to
     \Or(\Mm_{xz}^{u\#_\rho v}),\quad
     \rho\in[\rho_0,\infty).
$$
Here $\Mm_{xy}^u$
denotes the connected
component of $\Mm_{xy}$
containing $u$.
Let $[\dot u]$ denote
the orientation of $\Mm_{xy}^u$
provided by the flow.
The orientation
of a $k$-dimensional fibre
determined by an ordered $k$-tuple
of vectors is denoted by
$\langle v_1,\cdots ,v_k \rangle$.
Let
$[\langle v_1,\dots ,v_k \rangle]$
be the resulting
orientation of the whole
orientable vector bundle.
The map $\sigma^\#$ is defined
in case of flow orientations by
(see Figure~\ref{fig:fig-embd})
$$
     \sigma^\#\left([\dot u],[\dot v]\right)
     :=\left[\langle-\nabla f(p_\rho),
     -\tfrac{d}{d\rho} p_\rho\rangle\right],
$$
and in the general case by
\begin{equation}\label{eq:or-glu}
     \sigma^\#\left([\Mm_{xy}^u],
     [\Mm_{yz}^v]\right)
     :=ab\;\sigma^\#\left([\dot u],
     [\dot v]\right),
\end{equation}
where $a,b\in\{\pm1\}$
are determined by
$[\Mm_{xy}^u]=a[\dot u]$ and
$[\Mm_{yz}^v]=b[\dot v]$.

\begin{theorem}[Coherence]
\label{thm:coh-or}
The gluing map~(\ref{eq:or-glu})
and the orientations provided by
Proposition~\ref{prop:orientation}
are compatible in the sense that
$$
     \sigma^\#\left(
     [\Mm_{xy}^u]_{ind},
     [\Mm_{yz}^v]_{ind}\right)
     =[\Mm_{xz}^{u\#_\rho v}]_{ind}.
$$
\end{theorem}

\begin{proof}
Define $n_u\in\{\pm1\}$
by the identity
$[\Mm_{xy}^u]_{ind}=n_u[\dot u]$,
then
\begin{equation}\label{eq:coh-or}
\begin{split}
     \sigma^\#\left(
     [\Mm_{xy}^u]_{ind},
     [\Mm_{yz}^v]_{ind}\right)
    &=n_un_v\sigma^\#
     \left([\dot u],[\dot v]\right)\\
    &=n_un_v\left[\langle-\nabla f(p_\rho),
     -\tfrac{d}{d\rho} p_\rho\rangle\right].
\end{split}
\end{equation}
The second equality
holds by definition of
$\sigma^\#$.
To compare
the right hand side with
$[\Mm_{xz}^{u\#_\rho v}]_{ind}$
we need to relate
the induced orientations
of the bundles
$T\Mm_{xy}^u$, $T\Mm_{yz}^v$
and $T\Mm_{xz}^{u\#_\rho v}$.
Unfortunately, the base
manifolds do not have
a common point.
On the other hand, the point $y$
lies in the closure
of all three base manifolds
and all three tangent
bundles can be extended
to $y$. This is due to 
the existence of the
limits\footnote{
Here nondegeneracy of
$y$ is crucial and
$\dot u(+\infty)$
and $\dot v(-\infty)$
are eigenvectors
of the Hessian of $f$
at $y$ corresponding
to a positive
and a negative eigenvalue,
respectively.
}
(see~\cite[Lemma~B.5]{Schw93} and
\cite{KMP00})
$$
     \lim_{t\to\infty}
     \frac{\frac{d}{dt}\phi_tu}
     {\norm{\frac{d}{dt}\phi_tu}}
     =:\dot u(+\infty),\qquad
     \lim_{t\to-\infty}
     \frac{\frac{d}{dt}\phi_tv}
     {\norm{\frac{d}{dt}\phi_tv}}
     =:\dot v(-\infty).
$$
Repeatedly
using~(\ref{eq:splitting})
shows that the orientations
of the fibres over $y$
are related by
(we use the same
notation for the bundles
extended to $y$)
\begin{equation*}
\begin{split}
     [T_yW^u(x)]
    &=[T_y\Mm_{xy}^u]_{ind}
     \oplus
     [\Vv_y W^s(y)] \\
    &=[T_y\Mm_{xy}^u]_{ind}
     \oplus
     [T_y W^u(y)] \\
    &=[T_y\Mm_{xy}^u]_{ind}
     \oplus
     [T_y\Mm_{yz}^v]_{ind}
     \oplus
     [\Vv_y W^s(z)]
\end{split}
\end{equation*}
and (since 
$y=\lim_{\rho\to\infty}p_\rho$)
$$
     [T_yW^u(x)]
     =[T_y\Mm_{xz}^{u\#_\rho v}]_{ind}
     \oplus
     [\Vv_y W^s(z)].
$$
Hence
$$
     [T_y\Mm_{xz}^{u\#_\rho v}]_{ind}
     =[T_y\Mm_{xy}^u]_{ind}
     \oplus
     [T_y\Mm_{yz}^v]_{ind}
     =n_un_v [\dot u(+\infty)]
     \oplus[\dot v(-\infty)].
$$
The ordered pairs
$\langle \dot u(+\infty),
\dot v(-\infty)\rangle$
and $\langle-\nabla f(p_\rho),
-\tfrac{d}{d\rho} p_\rho\rangle$
represent the same
orientation of
$\Mm_{xz}^{u\#_\rho v}$.
By~(\ref{eq:coh-or})
this proves the theorem.
\end{proof}

\section{Morse homology}
\label{sec:Morse-homology}

\subsection{Morse-Witten complex}
\label{subsec:MW-complex}

\begin{definition}\rm
The \emph{Morse chain groups
associated to a
Morse function $f$,
with integer coefficients
and graded by the Morse index},
are the free abelian
groups generated by
the critical points of $f$
of Morse index $k$
$$
     \CM_k(M,f)
     := \bigoplus_{x\in\Crit_k f}
     \Z x,\qquad
     k\in\Z.
$$
A sum over
the empty set is
understood to be zero.
\end{definition}

The chain groups are finitely
generated by
Corollary~\ref{co:crit-finite}.
Let us choose a Riemannian
metric $g$ on $M$.
If $(g,f)$ is not a Morse-Smale pair,
replace it by a $C^1$-close
Morse-Smale pair according
to Theorem~\ref{thm:MS}.
Since the type of a Morse function
is locally constant, both chain groups
are canonically isomorphic.
From now on we assume
that $(g,f)$ is Morse-Smale.
Choose an
orientation for every
unstable manifold and
denote this set of choices
by $Or$.

\begin{definition}\rm
Assume $\IND_f(x)-\IND_f(y)=1$
and let $u\in\Hat\Mm_{xy}$.
The orbit $\Oo(u)$
is a connected component
of $\Mm_{xy}$
and hence carries the induced
orientation $[\Oo(u)]_{ind}$
provided by
Proposition~\ref{prop:orientation}.
Denoting the
flow orientation by
$[\dot u]$, the
\emph{characteristic sign}
$n_u=n_u(Or)$
is defined by
$$
     [\Oo(u)]_{ind}=n_u [\dot u].
$$
\end{definition}

\begin{definition}\rm
The \emph{Morse-Witten boundary operator}
$$
     \p_k=\p_k(M,f,g,Or):
     \CM_k(M,f)\to\CM_{k-1}(M,f)
$$
is given on a generator $x$ by
$$
     \p_k x
     :=\sum_{y\in \Crit_{k-1} f}\,
     n(x,y) y,\qquad
     n(x,y):=\sum_{u\in\Hat\Mm_{xy}}
     n_u,
$$
and extended to general chains
by linearity.
\end{definition}

\begin{figure}[ht]
\begin{minipage}[b]{.46\linewidth}
    \centering
    \epsfig{figure=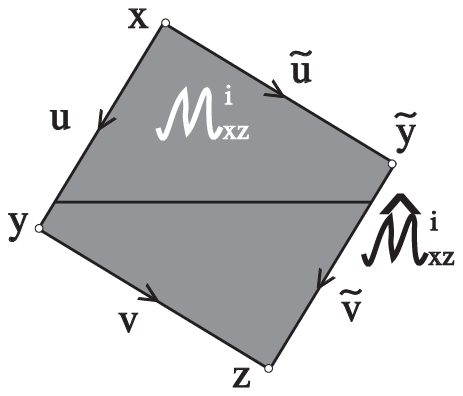}
    \caption{$\Hat\Mm_{xz}^i\simeq(0,1)$}
    \label{fig:fig-diam}
\end{minipage}
\hfill
\begin{minipage}[b]{.46\linewidth}
    \centering
    \epsfig{figure=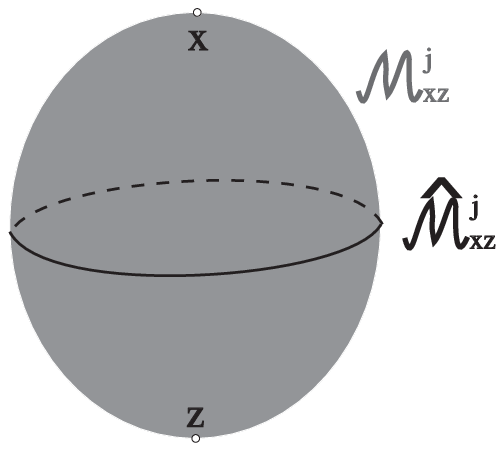}
    \caption{$\Hat\Mm_{xz}^j\simeq S^1$}
    \label{fig:fig-sphe}
\end{minipage}
\hfill
\end{figure}

Both sums in the definition of $\p_k$
are finite by
Corollary~\ref{co:crit-finite}
and
Proposition~\ref{pr:0D-compact},
respectively.
To prove that $\p$ satisfies
$\p^2=0$ we need to investigate
the 1-dimensional
components of the space
of connecting orbits.
Fix $x\in \Crit_k f$
and $z\in \Crit_{k-2} f$.
By Theorem~\ref{thm:manifold}
the orbit space
$\Hat\Mm_{xz}$ is a manifold
(without boundary)
of dimension 1 and therefore
its connected components $\Mm_{xz}^i$
are diffeomorphic either
to $(0,1)$ or to $S^1$
(see Figs.~\ref{fig:fig-diam}--\ref{fig:fig-sphe}).

\begin{proposition}
\label{pr:diamond}
Let $x\in \Crit_k f$
and $z\in \Crit_{k-2} f$,
then the following are true.
(i) The set 
\[
     \Bb^1_{xz}
     :=\{ (u,v)\mid
     \text{$u\in \Hat\Mm_{xy}$,
     $v\in \Hat\Mm_{yz}$,
     for some $y\in\Crit_{k-1} f$}
     \}
\]
of broken
orbits of order two
between $x$ and $z$
corresponds precisely
to the ends of the noncompact
connected components
of $\Hat\Mm_{xz}$.\\
(ii) Two broken
orbits $(u,v)$ and
$(\tilde{u},\tilde{v})$
corresponding
to the same connected component
$\Hat\Mm_{xz}^i$
are called \emph{cobordant}.
Their characteristic signs
satisfy
\[
     n_u n_v
     +n_{\tilde{u}} n_{\tilde{v}}
     =0.
\]
\end{proposition}

\begin{proof}
By Theorem~\ref{thm:compactness}
a connected component
$\Hat\Mm_{xz}^i\simeq(0,1)$
is compact up to
broken orbits of order two.
Hence we obtain two
broken orbits $(u,v)$
and $(\tilde{u},\tilde{v})$ --
one corresponding to each end.
The last statement in
Theorem~\ref{thm:gluing}
implies
$(u,v)\not=(\tilde{u},\tilde{v})$.
(However $u=\tilde{u}$ and
$v\not=\tilde{v}$, or vice versa,
is possible as
Example~\ref{ex:exp1} shows).
The Gluing Theorem~\ref{thm:gluing}
also tells that
each broken orbit $(u,v)$
corresponds to a noncompact
end of $\Hat\Mm_{xz}$.
This concludes the
proof of part (i).
To prove part (ii)
use the 
definition~(\ref{eq:or-glu})
of $\sigma^\#$,
Theorem~\ref{thm:coh-or},
and the
fact that $u\#_\rho v,
\tilde{u}\#_\rho \tilde{v}
\in\Hat\Mm_{xz}^i$
to obtain the following
identities for the
orientation of
$\Mm_{xz}^i$, 
namely
\begin{alignat*}{2}
     n_un_v[\langle\nabla f(p_\rho),
     \tfrac{d}{d\rho} p_\rho\rangle]
    &=n_un_v\sigma^\#\left([\dot u],
     [\dot v]\right)
   &&=\sigma^\#\left([u]_{ind},
     [v]_{ind}\right) \\
    &=[\Hat\Mm_{xz}^{u\#_\rho v}]_{ind}
   &&=[\Hat\Mm_{xz}^{\tilde{u}\#_\rho 
     \tilde{v}}]_{ind} \\
    &=\sigma^\#\left([\tilde{u}]_{ind},
     [\tilde{v}]_{ind}\right)
   &&=n_{\tilde{u}} n_{\tilde{v}}
     \sigma^\#\left([\dot {\tilde{u}}],
     [\dot {\tilde{v}}]\right) \\
    &=n_{\tilde{u}} n_{\tilde{v}}
     [\langle\nabla f(\tilde{p}_\rho),
     \tfrac{d}{d\rho} \tilde{p}_\rho\rangle]
    &&=-n_{\tilde{u}} n_{\tilde{v}}
     [\langle\nabla f(p_\rho),
     \tfrac{d}{d\rho} p_\rho\rangle].
\end{alignat*}
The last step follows, because
$\frac{d}{d\rho} \tilde{p}_\rho$
and $\frac{d}{d\rho} p_\rho$
both point outward
along the boundary
of $\Hat\Mm_{xz}^i$.
\end{proof}

\begin{theorem}[Boundary operator]
\label{thm:boundary}
$
     \p_{k-1} \p_k=0,\quad
     \forall k\in\Z.
$
\end{theorem}

\begin{proof}
By linearity,
definition of $\p$
and $\Bb^1_{xz}$, and
Proposition~\ref{pr:diamond}
it follows
\begin{equation*}
\begin{split}
     \p_{k-1} \p_k x
    &=\sum_{z\in \Crit_{k-2} f}
     \left(
     \sum_{y\in \Crit_{k-1} f}
     \sum_{u\in\Hat\Mm_{xy}}
     \sum_{v\in\Hat\Mm_{yz}}
     n_u n_v
     \right) z \\
    &=\sum_{z\in \Crit_{k-2} f}
     \left(
     \sum_{(u,v)\in\Bb^1_{xz}}
     n_u n_v
     \right) z \\
    &=\sum_{z\in \Crit_{k-2} f}
     \left(\sum_{
       \overset{\scriptstyle
       \text{connected components
       $\Hat\Mm_{xz}^i$}}
         {\text{of $\Hat\Mm_{xz}$
         diffeomorphic to $(0,1)$}}
     }\left(
     n_{u_i} n_{v_i}
     +n_{\tilde{u}_i} n_{\tilde{v}_i}
     \right)
     \right) z \\
    &=0.
\end{split}
\end{equation*}
\end{proof}

\begin{definition}\rm
Given a closed smooth
finite dimensional manifold $M$,
a Morse function $f$
and a Riemannian metric $g$ on $M$
such that the Morse-Smale condition holds,
denote by $Or$
a choice of orientations
of all unstable manifolds
associated to
the vector field $-\nabla f$.
Then the \emph{Morse homology groups
with integer coefficients}
are defined by
$$
     \HM_k(M;f,g,Or;\Z)
     :=\frac{\ker\p_k}
     {\im \p_{k+1}},\qquad
     k\in\Z.
$$
\end{definition}

In view of the Smale
Transversality
Theorem~\ref{thm:MS}
and the Continuation
Theorem~\ref{thm:continuation},
we can in fact define these homology groups
for every pair $(f,g)$, Morse-Smale or not,
and every choice of orientations $Or$.
Since they are all naturally isomorphic,
we shall denote them by
$\HM_*(M;\Z)$.

\begin{example}
\rm
Going back to
Example~\ref{ex:exp1}
let us calculate the characteristic
signs using the orientations
indicated in
Figure~\ref{fig:fig-exp1}.
We obtain
$$
     n_{u_1}=n_{u_2}
     =n_{\tilde{v}}=-1,\qquad
     n_v=+1.
$$
Hence
$\HM_2=\langle x_1-x_2 \rangle
\simeq\Z$,
$\HM_1=0$ and
$\HM_0=\langle z \rangle\simeq\Z$.
\end{example}

\subsection{Continuation}
\label{subsec:continuation}

In this section we present
Po\'zniak's~\cite{Po91}
construction of continuation
maps, i.e. of natural
grading preserving
isomorphisms
$$
     \Psi^{\beta\alpha}_*:
     \HM_*(M;f^\alpha,g^\alpha,Or^\alpha)
     \to
     \HM_*(M;f^\beta,g^\beta,Or^\beta)
$$
associated to any choice of
Morse-Smale pairs
$(f^\alpha,g^\alpha)$ and
$(f^\beta,g^\beta)$
and orientations
$Or^\alpha$ and $Or^\beta$
of all unstable manifolds.

\begin{theorem}[Continuation]
\label{thm:continuation}
$
     \Psi^{\beta\alpha}_*
     =(\Psi^{\alpha\beta}_*)^{-1},\quad
     \Psi^{\gamma\beta}_*\Psi^{\beta\alpha}_*
     =\Psi^{\gamma\alpha}_*.
$
\end{theorem}

The remaining part of this section
prepares the proof of the theorem.

\subsubsection*{The chain map
associated to a homotopy}

Let $(f^\alpha,g^\alpha,Or^\alpha)$
and $(f^\beta,g^\beta,Or^\beta)$
be as above.
Fix homotopies $\{f_s\}_{s\in[0,1]}$ and
$\{g_s\}_{s\in[0,1]}$ from
$f^\alpha$ to $f^\beta$ and
$g^\alpha$ to $g^\beta$, respectively.
By rescaling and smoothing the homotopies
near the endpoints, if necessary,
we may assume that they
are constant near the end points.
More precisely, for some fixed
$\delta\in[0,1/4]$ we have
$f_s\equiv f^\alpha$,
for $s\in[0,\delta]$, and
$f_s\equiv f^\beta$,
for $s\in[1-\delta,1]$, and similarly
for $g_s$. Such homotopies
$h^{\alpha\beta}
:=(f_s,g_s)$
are called \emph{admissible}.
Let us parametrize
$S^1$ by the interval $[-1,1]$
with endpoints being identified.

\begin{figure}[ht]
  \centering
  \epsfig{figure=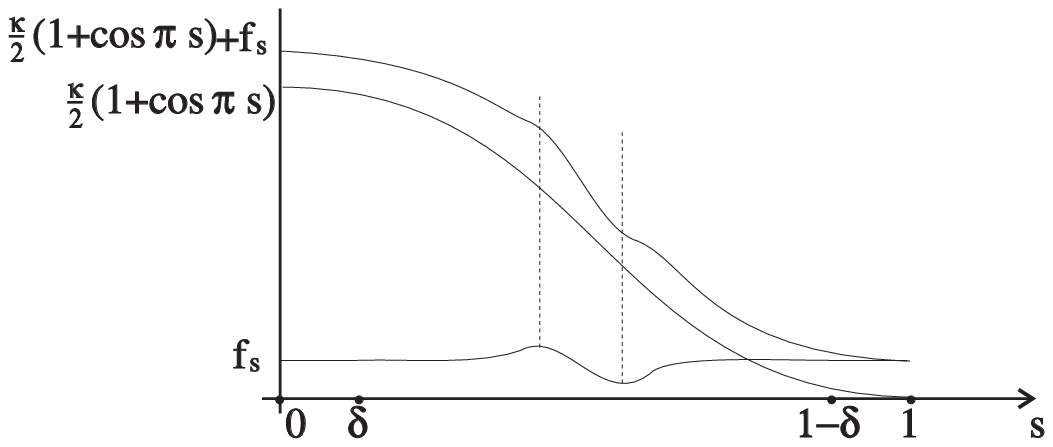}
  \caption{If $\kappa>>1$,
           then $\p_s F$ vanishes
           only at $s=0,1$
            }
  \label{fig:fig-mors}
\end{figure}

\begin{lemma}\label{le:F}
The set of critical points
of the function
$F=F_\kappa:M\times S^1\to\R$,
$$
     F(q,s)
     :=\frac{\kappa}{2}
     \left(1+\cos \pi s\right)
     +f_{\abs{s}}(q),
$$
coincides with
$\left(\Crit f^\alpha\times\{0\}\right)\cup
\left(\Crit f^\beta\times\{1\}\right)$,
for every positive real
\begin{equation}\label{eq:kappa}
     \kappa
     >
     \frac{2\max_{M\times S^1}\p_s f}
     {\pi\sin \pi(1-\delta)}.
\end{equation}
Then all critical points
are nondegenerate and their
Morse indices satisfy
$$
     \IND_F(x,0)=\IND_{f^\alpha}(x)+1,\qquad
     \IND_F(y,1)=\IND_{f^\beta}(y).
$$
\end{lemma}
The proof of the lemma is illustrated
by Figure~\ref{fig:fig-mors} and
left as an exercise.
Define a product metric
on $M\times S^1$ by
$$
     G_{q,s}:=g_{q,\abs{s}}\oplus 1
$$
and consider the negative
gradient flow of $F$ with
respect to $G$. We orient
all unstable manifolds by
taking the choice 
$$
     Or^{\alpha\beta}
     :=(Or^\alpha\oplus\p_s)
     \cup Or^\beta.
$$
The Morse-Witten complex
associated to
$(M\times S^1,F,G,Or^{\alpha\beta})$
has the following properties.
Its chain groups
split by Lemma~\ref{le:F}
\begin{equation}\label{eq:split}
     \CM_k(M\times S^1,F)
     \simeq
     \CM_{k-1}(M,f^\alpha) \oplus
     \CM_k(M,f^\beta).
\end{equation}
Setting $-\nabla F(q,s)
=:(\xi_{q,s},a_{q,s})
\in T_qM\times\R$
calculation shows
$$
     a_{q,s}
     =-\frac{\p}{\p s} F(q,s)
     =\frac{\kappa\pi}{2}\sin \pi s
     -\frac{\p}{\p s} f_{\abs{s}}(q)
$$
and $\xi_{q,s}$ is determined by the identity
$$
     df_{\abs{s}}(q) \cdot
     =
     g_{q,\abs{s}}(-\xi_{q,s},\cdot).
$$
By~(\ref{eq:kappa})
the zeroes of $a$
are precisely given by
the union of
$M\times\{0\}$
and $M\times\{1\}$.
In particular, both sets
are flow invariant and the restricted
flow coincides with the
one generated by
$-\nabla f^\alpha$ and
$-\nabla f^\beta$, respectively.
The characteristic
sign $n_u$ of an isolated
flow line $u\in\Hat\Mm_{xy}$
is related to the
characteristic signs
associated to the corresponding
flow lines $(u,0)$ and $(u,1)$
in $M\times\{0\}$ and $M\times\{1\}$,
respectively, by
$$
     n_{(u,1)}=n_u=-n_{(u,0)}.
$$
Replacing $\kappa$
by any number larger than
$(\max f^\beta-\min f^\alpha)$
guarantees that there are
no flow lines from
$M\times\{1\}$ to 
$M\times\{0\}$, because the
negative gradient flow
decreases along trajectories.

Concerning the Morse-Smale
condition we observe that
unstable and stable manifolds of
critical points $(x,i)$ and $(y,i)$,
respectively, intersect
transversally for $i=0$
and also for $i=1$.
However, this is not necessarily
the case for $W^u(x,0)$ and
$W^s(y,1)$. If necessary,
replace $(F,G)$
by a sufficiently
$C^1$-close Morse-Smale pair
(without changing notation).
The number and
Morse indices of critical points
as well as the
structure of connecting manifolds,
as long as they arise from
transversal intersections are preserved.
In particular, the two
chain subcomplexes sitting
at $M\times\{0\}$ and at
$M\times\{1\}$
remain unchanged.

Flow lines
from $M\times\{0\}$
to $M\times\{1\}$
converge at the ends
to $(x,0)$ and $(y,1)$,
where $x,y$ are
critical points of
$f^\alpha$ and $f^\beta$, respectively.
In case
$\IND_F(x,0)-\IND_F(y,1)=1$
there are finitely many
such flow lines.
The algebraic count of
those which pass through
$M\times\{1/2\}$
defines a map
$$
     \psi_*^{\beta\alpha}
     =\psi_*^{\beta\alpha}
     (h^{\alpha\beta},Or^{\alpha\beta}):
     \CM_*(M,f^\alpha)
     \to
     \CM_*(M,f^\beta)
$$
given on a generator $x$ by
\begin{equation}\label{eq:psi}
     x\mapsto
     \sum_{
       \overset{\scriptstyle
       (u,c)\in\Hat\Mm_{(x,0)(y,1)}}
         {\Oo((u,c))\cap (M\times\{1/2\})\not=
         \emptyset}
     }
     n_{(u,c)} y.
\end{equation}
(It is an exercise to
show that counting
all flow lines
produces the zero map).
Hence the boundary operator
$\Delta_k=\Delta_k(M\times S^1,
F,G,Or^{\alpha\beta})$
is of the form
$$
     \Delta_k
     \simeq
     \begin{pmatrix}
     -\p^\alpha_{k-1}&0\\
     \psi_{k-1}^{\beta\alpha}
     &\p^\beta_k
     \end{pmatrix}
$$
with respect to the
splitting~(\ref{eq:split}).
Theorem~\ref{thm:boundary}
states $\Delta_{k-1}\Delta_k=0$
and this implies
$
     \psi_{k-2}^{\beta\alpha}\p_{k-1}^\alpha
     =\p_{k-2}^\beta\psi_{k-1}^{\beta\alpha}
$.
Hence
$\psi_*^{\beta\alpha}(h^{\alpha\beta},Or^{\alpha\beta})$
is a chain map and we denote
the induced map on homology
by $[\psi_*^{\beta\alpha}(h^{\alpha\beta},
Or^{\alpha\beta})]$.
Observe that by $[\cdot]$
we also denote orientations.
However, the meaning
should be clear
from the context.

\begin{remark}[Constant homotopies]
\rm
Let $h^\alpha$ denote the
pair of constant homotopies
$(f^\alpha,g^\alpha)$
and set $Or^{\alpha\alpha}
=(Or^\alpha\oplus\p_s)
\cup Or^\alpha$.
Then the isolated flow lines
of the negative gradient flow
on $M\times S^1$ come in
quadruples and with characteristic
signs as shown
in Figure~\ref{fig:fig-quad}.
This implies
\begin{equation}\label{eq:const}
     \psi_*^{\alpha\alpha}(h^\alpha,Or^{\alpha\alpha})=\1.
\end{equation}
\begin{figure}[ht]
  \centering
  \epsfig{figure=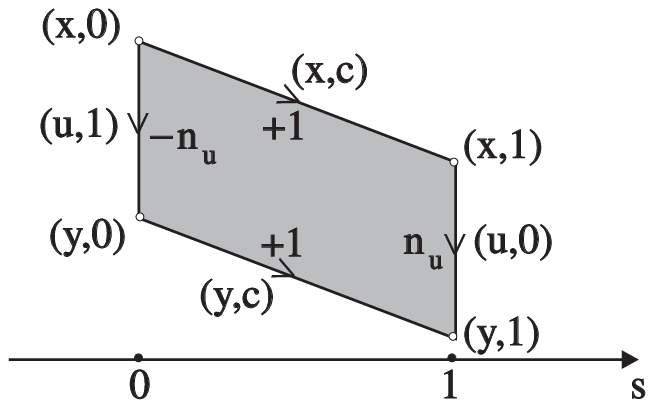}
  \caption{Isolated flow lines
           in case of constant homotopies}
  \label{fig:fig-quad}
\end{figure}
\end{remark}

\begin{remark}[Variation of metric]
\rm
Morse-Witten chain complexes
corresponding to Morse-Smale
pairs $(f,g^\alpha)$ and
$(f,g^\beta)$ are
\emph{chain} isomorphic,
as observed by Cornea and
Ranicki~\cite[Rmk.~1.23~(c)]{CR03}.
To see this reorder the generators
of $\CM_k(M,f)$ according
to descending value of $f$,
choose a homotopy $g_s$
and define $F$, $G$ and $Or^{\alpha\alpha}$
as above.
Now the key observation is that 
the chain homomorphism
$\psi_k^{\alpha\alpha}((f,g_s),Or^{\alpha\alpha})$
has diagonal entries $+1$
and only zeroes above the diagonal:
relevant isolated flow lines of
$-\nabla F$ from $M\times 0$
to $M\times 1$
are given by 
pairs $(u,c):\R\to M\times [0,1]$
satisfying
$$
     \begin{pmatrix}
     \dot u\\\dot c
     \end{pmatrix}
     =
     \begin{pmatrix}
     \xi_{u,c}\\
     \tfrac{\kappa\pi}{2}
     \sin \pi c
     \end{pmatrix},\quad
     df|_u(\cdot)
     =g_{u,\abs{c}}(-\xi_{u,c},\cdot),\quad
     \lim_{t\to\pm\infty}(u,c)
     \in \Crit_k f\times \{0\}.
$$
For every $x\in\Crit_k f$
there is a unique
isolated flow line
connecting $(x,0)$ to $(x,1)$.
It is of the form
$(x,c)$ and transversality is
automatically satisfied
(therefore it survives small
perturbations possibly required to
achieve Morse-Smale transversality
for $(F,G)$ and actually
defining $\psi_k^{\alpha\alpha}$).
Hence every
diagonal entry of $\psi_k^{\alpha\alpha}$
equals $+1$.
The identity
for flow lines $(u,c)$
$$
     \frac{d}{dt}f(u(t))
     =df|_{u(t)}\dot u(t)
     =g_{u(t),\abs{c(t)}}
      (-\xi_{u(t),c(t)},\dot u(t))
     =-\Abs{\dot u(t)}_{g_{u(t),\abs{c(t)}}}^2
$$
shows that $f$ is strictly decreasing
along nonconstant $u$-components,
hence all elements of
$\psi_k^{\alpha\alpha}$ strictly 
above the diagonal are zero.
\end{remark}

\subsubsection*{Homotopies of homotopies}

Given four Morse functions,
metrics and choices of orientations
$(f^i,g^i,Or^i)$,
$i=\alpha,\beta,\gamma,\delta$,
and homotopies of homotopies
$\{f_{s,r}\}_{s,r\in[0,1]}$ and
$\{g_{s,r}\}_{s,r\in[0,1]}$
satisfying, for some fixed $\eps\in[0,1/4]$,
$$
     f_{s,r}
     =\begin{cases}
     f^\alpha &, (s,r)\in[0,\eps]\times[0,\eps] \\
     f^\beta  &, (s,r)\in[1-\eps,1]\times[0,\eps] \\
     f^\gamma &, (s,r)\in[0,\eps]  \times[1-\eps,1] \\
     f^\delta &, (s,r)\in[1-\eps,1]\times[1-\eps,1]
     \end{cases}
$$
and similarly for $g_{s,r}$,
we define
a function and metric on
$M\times S^1\times S^1$ by
\begin{equation*}
\begin{split}
     F(q,s,r)
    &:=\frac{\kappa}{2}
     \left(1+\cos \pi s\right)
     +\frac{\rho}{2}
     \left(1+\cos \pi r\right)
     +f_{\abs{s},\abs{r}}(q)\\
     G_{q,s,r}
    &:=g_{q,\abs{s},\abs{r}}\oplus1\oplus1.
\end{split}
\end{equation*}
Choosing $\kappa,\rho>0$ sufficiently
large it follows as in
Lemma~\ref{le:F} that $F$ is
Morse and there is a splitting
\begin{equation}\label{split2}
\begin{split}
    &\CM_k(M\times S^1\times S^1,F)\\
    &\simeq
     \CM_{k-2}(M,f^\alpha) \oplus
     \CM_{k-1}(M,f^\beta)
     \oplus \CM_{k-1}(M,f^\gamma) 
     \oplus \CM_k(M,f^\delta).
\end{split}
\end{equation}
Orient all unstable manifolds
of $-\nabla F$ on
$M\times S^1\times S^1$ by
$$
     Or^{\alpha\beta\gamma\delta}
     :=(Or^\alpha\oplus\p_s\oplus\p_r)
     \cup (Or^\beta\oplus\p_r)
     \cup (Or^\gamma\oplus\p_s)
     \cup Or^\delta.
$$
Arguing similarly as in the
former subsection we conclude
that the boundary operator
is, with respect to the
splitting~(\ref{split2}),
represented by
$$
     \Delta_k
     (M\times S^1\times S^1,
     F,G,Or^{\alpha\beta\gamma\delta})
     \simeq
     \begin{pmatrix}
     \p^\alpha_{k-2}&0&0&0\\
     \psi^{\beta\alpha}_{k-2}
     &-\p^\beta_{k-1}&0&0\\
     -\psi^{\gamma\alpha}_{k-2}
     &0&-\p^\gamma_{k-1}&0\\
     \Lambda^{\delta\alpha}_{k-2}
     &\psi^{\delta\beta}_{k-1}
     &\psi^{\delta\gamma}_{k-1}
     &\p^\delta_k
     \end{pmatrix},
$$
where $\Lambda$ and the $\psi$'s
are defined similar to~(\ref{eq:psi})
by counting isolated flow lines.
The signs arise as follows.
The characteristic sign
$n_u$ of an isolated flow line $u$ in
one of the subcomplexes, e.g.
the first summand
in~(\ref{split2}),
and the sign $n_{(u,0,0)}$ 
of the corresponding
element of the full complex
are related by
$$
     n_{(u,0,0)}=n_u^\alpha,\quad
     n_{(u,1,0)}=-n_u^\beta,\quad
     n_{(u,0,1)}=-n_u^\gamma,\quad
     n_{(u,1,1)}=n_u^\delta.
$$
Similarly the signs $n_{(u,c)}$
in the definition~(\ref{eq:psi}) of
$\psi_*$ and the corresponding ones
in the full complex
are related by
$$
     n_{(u,c,0)}=n_{(u,c)}^{\beta\alpha},\quad
     n_{(u,c,1)}=n_{(u,c)}^{\delta\gamma},\quad
     n_{(u,0,c)}=-n_{(u,c)}^{\gamma\alpha},\quad
     n_{(u,1,c)}=n_{(u,c)}^{\delta\beta}.
$$
The identity
$
     \Lambda^{\delta\alpha}_{k-3}
     \p^\alpha_{k-2}
     +\p^\delta_{k-1}
     \Lambda^{\delta\alpha}_{k-2}
     =\psi^{\delta\gamma}_{k-2}
     \psi^{\gamma\alpha}_{k-2}
     -\psi^{\delta\beta}_{k-2}
     \psi^{\beta\alpha}_{k-2}
$
is a consequence of
$\Delta_{k-1}\Delta_k=0$
and shows that $\Lambda$
is a chain homotopy between
$\psi^{\delta\gamma}_{k-2}
\psi^{\gamma\alpha}_{k-2}$
and
$\psi^{\delta\beta}_{k-2}
\psi^{\beta\alpha}_{k-2}$.
Therefore
\begin{equation}\label{eq:homhom}
     \bigl[\psi^{\delta\gamma}_*\bigr]
     \bigl[\psi^{\gamma\alpha}_*\bigr]
     =
     \bigl[\psi^{\delta\beta}_*\bigr]
     \bigl[\psi^{\beta\alpha}_*\bigr].
\end{equation}
It is important to recall
that $\psi^{\beta\alpha}_*$
actually abbreviates
$\psi_*^{\beta\alpha}(h^{\alpha\beta},Or^{\alpha\beta})$.

\begin{proposition}\label{pr:indep}
The induced map on homology
$$
     \Psi^{\beta\alpha}_*
     :=\bigl[\psi_*^{\beta\alpha}(h^{\alpha\beta},
     Or^{\alpha\beta})\bigr]:
     \HM_*(M;f^\alpha,g^\alpha,Or^\alpha)
     \to
     \HM_*(M;f^\beta,g^\beta,Or^\beta)
$$
is independent of the choice of
$h^{\alpha\beta}$
and
$Or^{\alpha\beta}$
and satisfies
$\Psi^{\alpha\alpha}_*=\1$.
\end{proposition}

\begin{figure}[ht]
\begin{minipage}[b]{.48\linewidth}
    \centering
    \epsfig{figure=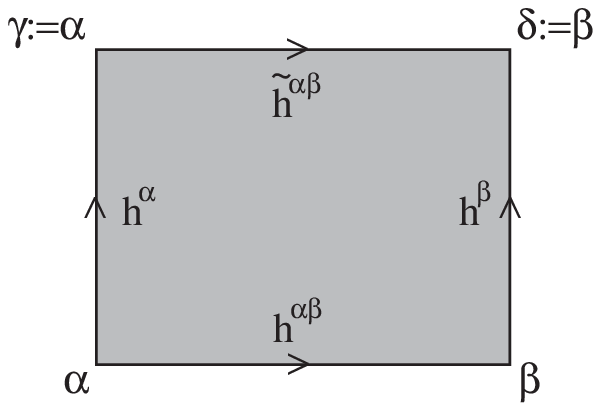}
    \caption{Proof of
            Proposition~\ref{pr:indep}}
    \label{fig:fig-hom1}
\end{minipage}
\hfill
\begin{minipage}[b]{.46\linewidth}
    \centering
    \epsfig{figure=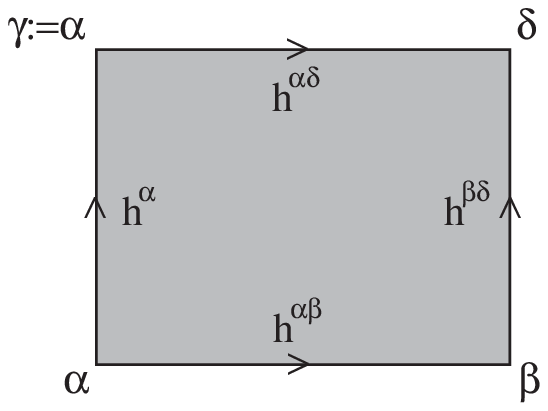}
    \caption{Proof of
             Theorem~\ref{thm:continuation}}
    \label{fig:fig-hom2}
\end{minipage}
\hfill
\end{figure}

\begin{proof}
In~(\ref{eq:homhom})
choose $\gamma:=\alpha$,
$\delta:=\beta$ and homotopies
as indicated in Figure~\ref{fig:fig-hom1},
where
$h^\alpha=(f^\alpha,g^\alpha)$
denotes constant homotopies.
By~(\ref{eq:homhom}) we have
$$
     \bigl[\psi_*^{\beta\alpha}
       (\tilde{h}^{\alpha\beta},Or^{\alpha\beta})
    \bigr]
     \bigl[\psi_*^{\alpha\alpha}
       (h^{\alpha},Or^{\alpha\alpha})
     \bigr]
     =
     \bigl[\psi_*^{\beta\beta}
       (h^{\beta},Or^{\beta\beta})
     \bigr]
     \bigl[\psi_*^{\beta\alpha}
       (h^{\alpha\beta},Or^{\alpha\beta})
    \bigr],
$$
and~(\ref{eq:const}) concludes the proof.
\end{proof}

\begin{proof}
[Proof of Theorem~\ref{thm:continuation}]
Choosing $\gamma:=\alpha$
in~(\ref{eq:homhom})
and homotopies as indicated
in Figure~\ref{fig:fig-hom2} proves
$\Psi_*^{\delta\beta}
\Psi_*^{\beta\alpha}
=\Psi_*^{\delta\alpha}$.
Setting $\delta:=\alpha$
concludes the proof.
\end{proof}

\subsection{Computation}
\label{subsec:computation}

Motivated by our examples
in Figures~\ref{fig:fig-exp1},
\ref{fig:fig-tor2} and~\ref{fig:fig-sphe}
one might conjecture

\begin{theorem}\label{thm:computation}
$\HM_*(M;\Z)\simeq \Ho_*^{sing}(M;\Z)$.
\end{theorem}

The examples also suggest
that the unstable manifolds
corresponding to a Morse cycle
are closely related to singular cycles.
This observation
is the origin of various
geometric approaches
towards a proof of the theorem.

However, the first proofs given
were somewhat less
geometrical; see Milnor~\cite[Lemma~7.2,
Theorem~7.4]{Mi65}
in case of self-indexing $f$;
for the general case see
Floer~\cite{Fl89},
Salamon~\cite{Sa90} and
Schwarz~\cite{Schw93}.

The geometric idea above
was made precise by
Schwarz~\cite{Schw99}
constructing pseudo-cycles
in an intermediate step.

Another idea for
a geometric proof has been
around for many years,
but to the best of my knowledge
never made rigorous:
a \emph{smooth triangulation of $M$}
is a pair $(K,h)$ of a simplicial complex
$K$ and a homeomorphism $h$
between the union $\abs{K}$ of all simplices
of $K$ and $M$, such that $h$ restricted
to any top-dimensional simplex is
an embedding.
It exists by~\cite{Ca35,Wh40}.
Given $(K,h)$,
the idea is to construct
a Morse function $f$ whose critical
points are precisely at the
barycenters\footnote{For
simplicity we shall
call the images of barycenters
under $h$ again barycenters.}
of the simplices and
a metric $g$,
such that $(f,g)$ is a Morse-Smale pair and
its negative gradient flow
has the following property.
From each barycenter of
a $k$-face there is precisely
one isolated flow line to 
the barycenter
of every $(k-1)$-face
in the boundary
of the $k$-face. The orientations
of the unstable manifolds
have to be chosen, such that
the characteristic signs
reflect the signs in the
definition of the simplicial
boundary operator.
Then the simplicial
and Morse chain complexes are identical
and the theorem follows
immediately.
A difficulty
in this approach is the fact
that top-dimensional simplices
do not in general fit
together smoothly.

Yet another idea is to
give $M$ the structure
of a CW-space viewing the
unstable manifolds
of a self-indexing Morse function
as cells and then relate
Morse homology to cellular homology.
In case that the metric is euclidean
with respect to the local
coordinates provided near critical
points by the Morse Lemma
this has been achieved
by Laudenbach~\cite{Ld92}.
In the general case
the difficulty arises
how to extend the
natural identification
of an open disc with
an unstable manifold,
which is provided by the flow,
\emph{continuously}
to the boundary:
some point might
converge to one critical
point, but arbitrarily close points
to another one.

\subsection{Remarks}
\label{subsec:remarks}

\subsubsection*{Morse inequalities}

To prove the 
Morse-inequalities~(\ref{eq:M-ineq})
we may assume without
loss of generality that
the homology of $M$ is 
torsion free (otherwise
choose rational coefficients
and observe that
$b_k(M;\Z)=b_k(M;\Q)$;
see e.g.~\cite[Satz 10.6.6]{SZ88}).
The dimension
of the free module
$\CM_k(M,f)$
was denoted by $c_k$.
Since the Morse-Witten boundary
operator $\p_k$
is a module homomorphism,
there is a splitting of
$
     \CM_k(M,f)
$
into the sum
of free submodules
$
     \ker \p_k\oplus
     \im \p_k
$.
This implies
\begin{equation}\label{eq:sum}
     c_k=\gamma_k+\beta_k,\qquad
     \gamma_k:=\dim \ker \p_k,\quad
     \beta_k:=\dim \im \p_k,
\end{equation}
with $\beta_0=0=\beta_{n+1}$.
Moreover, by Theorem~\ref{thm:computation}
\begin{equation}\label{eq:ineq}
     b_k
     =\rank \left(
     \frac{\ker\p_k}{\im \p_{k+1}}
     \right)
     = \gamma_k-\beta_k.
\end{equation}
The last identity is due to the
assumption that there is no torsion.
It follows by induction
using~(\ref{eq:sum})
and~(\ref{eq:ineq}) that
$$
     c_k-c_{k-1}+\dots\pm c_0
     =\beta_{k+1}+ b_k-b_{k-1}+\dots\pm b_0,\qquad
     k=0,\dots,n,
$$
and this proves~(\ref{eq:M-ineq}).

\subsubsection*{Relative Morse homology}

Given a Morse function $f$ with
regular values $a$ and $b$, one can
modify the definition of the
Morse chain groups
by using as generators only those
critical points $x$ with
$f(x)\in[a,b]$.
The resulting Morse homology
groups represent relative
singular homology. More
precisely, with
$M^a_f:=\{f\le a\}$ it holds
$$
     \HM_*^{(a,b)}(M,f,g,Or)
     \simeq \Ho_*^{sing}
     (M^b_f,M^a_f;\Z).
$$
Independence of $f$
holds true for such $f$
which are connected
by a \emph{monotone}
homotopy $f_s$
(meaning that $\p_s f_s\le 0$
pointwise).

Another point of view
is as follows.
Given a manifold $M$ with boundary
and a negative gradient
vector field $-\nabla f$
which is nonzero
on $\p M$ and points outward
on some boundary components,
denoted by $(\p M)^{out}$,
then
$$
     \HM_*(M,f,g,Or)
     \simeq \Ho_*^{sing}
     (M,(\p M)^{out};\Z).
$$

As an example consider
a 2-sphere in $\R^3$
as in Figure~\ref{fig:fig-sphe},
equipped with the height
function $f$ and the metric
induced by the euclidean metric
on $\R^3$. Denote the upper
hemisphere by $D_u$
and the lower one by $D_\ell$.
It is then easy to calculate
the absolute and relative
homology groups
of a disk $D$, namely
\begin{equation*}
\begin{split}
     \Ho_0^{sing}(D;\Z)
     \simeq \HM_0(D_\ell,-\nabla f)
    &=\langle z \rangle =\Z, \\
     \Ho_2^{sing}(D,\p D;\Z)
     \simeq \HM_2(D_u,-\nabla f)
    &=\langle x \rangle =\Z.
\end{split}
\end{equation*}
For all other values of $k$
these homology groups are zero,
since there are simply no generators
on the chain level.

\subsubsection*{Morse cohomology}

Let $(f,g)$ be a Morse-Smale pair.
The Morse chain groups $\CM_k(M,f)$
and cochain groups
$\CM^k(M,f):=\Hom(\CM_k(M,f),\Z)$
can be identified via
the map $\Crit_k f\ni x\mapsto\eta_x$,
where $\eta_x(z)$ is one in
case $z=x$ and zero otherwise.
Define the Morse
coboundary operator
$$
     \delta^k=\delta^k(M,f,g,Or)
     :\CM^k(M,f)\to\CM^{k+1}(M,f)
$$
on a generator $y$ by
counting isolated flow
lines of the negative
gradient flow ending
at $y$ (or equivalently
those of the \emph{positive}
gradient flow emanating
from $y$)
$$
     \delta^k y
     :=\sum_{x\in\Crit_{k+1} f} n(x,y) x.
$$
The homology of this
chain complex is called
\emph{Morse cohomology}
and denoted by
$\HM^*(M,f,g,Or)$.
It is naturally isomorphic
to $\Ho^*_{sing}(M;\Z)$.

\subsubsection*{Poincar\'e duality}

Consider the natural
identifications
\begin{equation}\label{eq:coho}
     \CM^k(M,f)
     \simeq \CM_k(M,f)
     \simeq \CM_{n-k}(M,-f).
\end{equation}
Assume that
$M$ is orientable
and fix an orientation
(or restrict to
coefficients
in $\Z_2$).
Taking a choice $Or$
of orientations of
all unstable manifolds
with respect to $-\nabla f$
determines orientations
of the stable manifolds,
which equal the unstable
manifolds with respect to
$-\nabla (-f)$. Let
$Or^{-f}=Or^{-f}(Or,[M])$
denote the
induced orientations.
Under the
identification~(\ref{eq:coho})
the boundary operators
$\delta^k(M,f,g,Or)$ and
$\p_{n-k}(M,-f,g,Or^{-f})$
coincide. Hence
$$
     \HM^k(M,f,g,Or)
     \simeq
     \HM_{n-k}(M,-f,g,Or^{-f})
     \simeq
     \HM_{n-k}(M,f,g,Or).
$$

\subsection{Real projectice space}
\label{subsec:real-proj-space}

\begin{figure}[ht]
  \centering
  \epsfig{figure=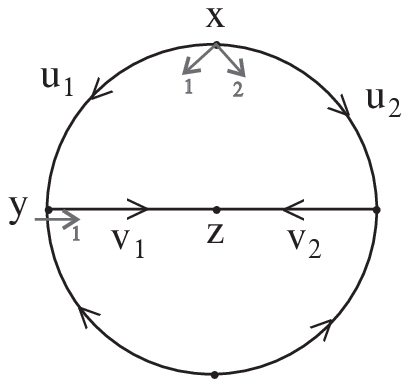}
  \caption{Morse function on $\RP^2$}
  \label{fig:fig-RP-2}
\end{figure}

Real projective
space exhibits enough
asymmetry in its (co)homology
to be a good testing ground
for the (rather informal)
remarks above.
View $\RP^2$ as the unit
disc in $\R^2$ with
opposite boundary points
identified and
consider a Morse function
as in
Figure~\ref{fig:fig-RP-2}
(cf.~\cite[Ch.~6, \S~3, Exc.~6]{Hi76})
having precisely three critical
points $x,y,z$ of Morse indices 2,1,0, 
respectively. 

\noindent
Hence
$$
     \CM_2=\CM^2=\langle x\rangle,\quad
     \CM_1=\CM^1=\langle y\rangle,\quad
     \CM_0=\CM^0=\langle z\rangle.
$$
Let the orientations
of the unstable manifolds
be as in
Figure~\ref{fig:fig-RP-2},
then
$$
     n_{u_1}=n_{u_2}=n_{v_1}=+1,\qquad
     n_{v_2}=-1,
$$
and the (co)boundary operators
$\p_k$ and $\delta^k$ act by
\begin{alignat*}{2}
     \p_2 x
    &=2y,\qquad\quad
    &\delta^2 x
    &=0,
     \\
     \p_1 y
    &=0,\qquad\quad
    &\delta^1 y
    &=2x,
    \\
     \p_0z
    &=0,\qquad\quad
    &\delta^0 z
    &=0.
\end{alignat*}
Hence integral
Morse (co)homology is given by
\begin{alignat}{4}\label{eq:rp-2}
     \HM_2(\RP^2;\Z) \nonumber
    &=0,\quad\qquad
    &\HM^2(\RP^2;\Z)
    &=\Z_2,\quad
     \\
     \HM_1(\RP^2;\Z)
    &=\Z_2,\quad\qquad
    &\HM^1(\RP^2;\Z)
    &=0,\quad
    \\
     \HM_0(\RP^2;\Z) \nonumber
    &=\Z,\quad
    &\HM^0(\RP^2;\Z)
    &=\Z,
\end{alignat}
whereas all (co)homology
groups with $\Z_2$-coefficients
equal $\Z_2$.

It is also interesting
to check in this example
that, if we replaced
in the Morse complex setup
the negative
by the positive
gradient flow (thereby
obtaining an ascending
boundary operator),
then the homology
of this new chain complex
would reproduce
the cohomology groups
in~(\ref{eq:rp-2}).

It is an instructive exercise
to extend the Morse function
in Figure~\ref{fig:fig-RP-2}
to $\RP^3$, viewed as the unit
three disc in $\R^3$ with
opposite boundary points identified,
and calculate its Morse (co)homology.

\end{document}